\documentclass{amsart}
\usepackage{psfig,amsmath,amssymb,latexsym}

\def\proof{\medskip\noindent{\sc proof. }}
\def\proofo{\medskip\noindent{\bf Proof of Theorem 7.1. }}

\def\EOP{\hfill$\Box$}
\def\natnum{\hbox{\rm I\kern-.17em N}}
\def\integ{\hbox{\rm Z\kern-.3em Z}}
\def\Ker{{\rm Ker}}
\def\Image{{\rm Im}}
\def\Tor{{\rm Tor}}
\def\Coker{{\rm Coker}}
\def\gr{{\rm gr}}
\def\deg{{\rm deg}}

\def\init{{\rm in}}
\def\frl{{\mathfrak L}}
\def\umm{{\mathbf m}}
\def\xx{{\bf \underline{x}}}
\def\field{{\mathbb K}}

\newtheorem{thm}{Theorem}[section]
\newtheorem{lem}[thm]{Lemma}
\newtheorem{prop}[thm]{Proposition}
\newtheorem{defn}[thm]{Definition}
\newtheorem{cor}[thm]{Corollary}
\newtheorem{qn}[thm]{Question}
\newtheorem{examp}[thm]{Example}
\newtheorem{obs}[thm]{Observation}
\newtheorem{rk}[thm]{Remark}

\begin{document}

\title[Gr\"obner Basis Bounds and Discrete Morse Theory]{Gr\"obner basis degree 
bounds on $\Tor^{k[\Lambda ]}_\bullet(k,k)_\bullet$
and Discrete Morse Theory for Posets}

\author[Hersh]{Patricia Hersh}
\author[Welker]{Volkmar Welker}
\date{}
\thanks{The first author was supported by an NSF Math Sciences Postdoctoral 
Research Fellowship. The second other was supported by
EU Research Training Network "Algebraic Combinatorics in Europe, " grant
HPRN-CT-2001-00272.}

\begin{abstract} 
The purpose of this paper is twofold.  
 
\noindent $\triangleright$ We give combinatorial 
bounds on the ranks of the groups 
$\Tor^{R}_\bullet(k,k)_\bullet$ 
in the case where $R = k[\Lambda]$ is an  affine semi-group ring,  
and in the process provide combinatorial proofs for bounds  
by Eisenbud, Reeves and Totaro on which Tor groups vanish. 
In addition, we show  
that if the bounds hold for a field $k$ then they hold for $\field[\Lambda]$ 
and any field $\field$.   
Moreover, we provide a combinatorial construction for a free resolution of  
$\field$ over $\field[\Lambda]$ which achieves these bounds. 
 
\noindent $\triangleright$ We extend the lexicographic discrete Morse  
function construction of Babson and Hersh for the determination of the  
homotopy type and homology of order complexes 
of posets to a larger class of facet orderings that includes orders  
induced by monomial term orders.  
 
Since it is known that the order complexes of finite intervals in the  
poset of monomials in $k[\Lambda]$ 
ordered by divisibility in $k[\Lambda ]$  
govern the $\Tor$-groups, the newly  
developed tools are applicable and serve as the main    
ingredients for the proof of the bounds and the construction of  
the resolution. 
\end{abstract} 
 
\maketitle 
 
\section{Introduction.} 
 
Let $\Lambda$ be a submonoid of $\natnum^e$ which is finitely generated by  
$\alpha_1, \ldots, \alpha_n$, 
and denote by $k[\Lambda ]$ the affine semi-group 
ring of $\Lambda$ generated over the field $k$ by  
monomials $\xx^{\alpha_i} = x_1^{\alpha_{i1}} 
\cdots x_e^{\alpha_{ie}}$, $ 1 \leq i \leq n$.  
Thus, $k[\Lambda] \cong k[z_1,\dots ,z_n]/I_{\Lambda }$ is the  
coordinate ring  
of an affine, not necessarily normal, toric variety. 
The isomorphism results from sending $z_i$ to $\xx^{\alpha_i}$,  
$1 \leq i \leq n$, and letting the toric ideal $I_{\Lambda }$ record the  
syzygies among the generators.  
The monoid $\Lambda $ is endowed with a partial order given by  
$\mu \leq \lambda $ if and only if $\lambda - \mu \in \Lambda $.  Denote 
by $\Delta (\mu ,\lambda )$ the simplicial complex of linearly ordered 
subsets $\mu < \mu_0 < \cdots < \mu_i < \lambda$ of the interval $[\mu,\lambda] 
:= \{ \gamma \in \Lambda~|~\mu \leq \gamma \leq \lambda\}$. 
 
Based on work by Laudal and Sletj\o e \cite{LS} and  Peeva, Reiner,  
Sturmfels [PRS], several recent papers (see e.g. [HRW], [BW])  
have used this partial order on $\Lambda$ 
and the simplicial homology of  
order complexes $\Delta(\mu,\lambda)$ 
for $\mu \leq \lambda$ in $\Lambda$ 
as a tool for understanding minimal free resolutions of the field $k$  
as a $k[\Lambda]$-module.   
In general, the minimal free resolution of the field 
$k$ over a $k$-algebra $R$ is still a mysterious object (see \cite{Av}) and  
even results known to hold by algebraic arguments pose hard and interesting  
combinatorial questions 
when $R = k[\Lambda]$ (see \cite{BjWe}).  Notably the Koszul property has  
attracted a lot of interest. This property is equivalent  
by work of Peeva, Reiner, Sturmfels [PRS] to 
the property of all intervals in the poset $\Lambda$ being 
Cohen-Macaulay over $k$. 
In general, it is known that a standard graded $k$-algebra is  
Koszul whenever its defining 
ideal has a quadratic Gr\"obner basis. In Peeva, Reiner, Sturmfels [PRS] 
and subsequent work \cite{HRW} a combinatorial  
understanding of this implication is 
developed.  In [BW] it is shown that if each interval in $\Lambda$ is  
shellable then 
it is actually possible to construct a minimal free resolution for $k$ as a 
$k[\Lambda ]$-module.   
 
In Section \ref{cm-section} we give an alternative 
combinatorial approach based on a discrete Morse function 
that also explains all these phenomena related to  
the Koszulness of $k[\Lambda]$ without requiring a shelling.  The main  
idea behind our combinatorial approach is quite natural, and is  
explained in Remark  
~\ref{intuitive-idea-remark} and the discussion that follows, after  
suitable notation is introduced.   
Section ~\ref{section-quad-applic} uses the discrete Morse function of 
Section ~\ref{cm-section} to  
provide a minimal free resolution for $k$ as a $k[\Lambda ]$-module 
when $I_{\Lambda }$ has a quadratic Gr\"obner basis, 
whether or not each interval in $\Lambda $ is shellable. 
 
In Section \ref{high_deg} we give the proof of our main result, a discrete 
Morse function on the order complex of $\Lambda $ 
which provides combinatorial upper bounds on all the multigraded Tor groups 
for $k[\Lambda ]$, most notably yielding the following: 
 
\begin{thm}\label{betti-vanish} 
Let $\Lambda \subseteq \natnum^e$ be an affine semi-group 
generated by $n$ elements of $\natnum^e$. 
Assume there is a field $k$ such that for $k[\Lambda] \cong 
k[z_1,\ldots, z_n]/I_\Lambda$ the ideal $I_\Lambda$ has a  
Gr\"obner basis of degree $d$, then 
\begin{itemize} 
\item[(i)] $\tilde{H}_i(\Delta(\hat{0},\lambda);\field )=0$ for  
$i< -1 + \frac{\deg (\lambda )-1}{d-1}$ and any field $\field$. 
\item[(ii)] $\Tor_i^{k[\Lambda ]}(\field,\field)_{\lambda }=0$  
for $i< 1+ \frac{\deg (\lambda )-1}{d-1}$ and any field $\field$. 
\end{itemize} 
Moreover, the vanishing of $\Tor$-groups  
is achieved by a free cellular resolution resulting from a discrete Morse 
function on $\Lambda$. 
\end{thm} 
 
Here we denote by $\deg(\lambda)$ the length, i.e. cardinality minus one,   
of a saturated chain in the poset interval $[\hat{0},\lambda]$. 
Note that if all generators of $\Lambda$ lie on an affine hyperplane then  
this grading actually makes $\Lambda$ a graded poset and $k[\Lambda]$ a  
standard graded $k$-algebra. In general, $\deg(\lambda)$ is the degree of  
the image of $\xx^\lambda$ in the associated graded ring of $R$.  
 
By results of [LS], Theorem \ref{betti-vanish} (i) will immediately 
imply Theorem \ref{betti-vanish} (ii). 
The vanishing of $\Tor$-groups in the case $\field = k$ also follows for  
general standard graded $k$-algebras, by a flat 
degeneration argument, from a result of Eisenbud, Reeves and Totaro  
(cf. [ERT]) about monomial ideals.  Our arguments are completely  
combinatorial.  
 
The main tool for the proof of Theorem \ref{betti-vanish}  
is discrete Morse theory, which was developed in the mid 90's by
Forman [Fo]. Discrete Morse theory is a tool
for determining the homology and homotopy type of a simplicial complex,
or more general a regular CW-complex.
In [BH] the authors develop tools that facilitate the use of 
discrete Morse theory in the case when the simplicial complex is the order
complex of a poset whose edges are labeled.  
The edges $\mu < \lambda$ of the Hasse diagram  
of $\Lambda$ are naturally labeled by $\lambda - \mu$ which by the definition 
of the order relation is one of the generators $\lambda_1, \ldots, 
\lambda_n$ of $\Lambda$. Thus, to 
any saturated chain $\mu_0 < \cdots < \mu_i$ in $\lambda$ there is 
associated the monomial $z_{j_1} \cdots z_{j_n}$ in  
$k[z_1,\ldots, z_n]$, where $\mu_i - \mu_{i-1} = \lambda_{j_i}$. 
In particular, any term order on $k[z_1, \ldots , z_n]$ induces a partial order 
on the finite saturated chains in $\Lambda$.
In the case of Gr\"obner bases with properties analogous to requirements made in
 [PRS] (called `supporting a poset''), we can apply the results from [BH]
on lexicographic discrete Morse function which then  
easily give a degree $d$ analogue of a lexicographic shelling, and imply the  
desired connectivity bound. In this Morse function, collections of at 
most $d$ labels, given by descents and Gr\"obner basis leading terms,  
will play the role traditionally filled by the descents in a  
lexicographic shelling.  In order to able to handle arbitrary Gr\"obner bases, 
whether or not they support a  poset, we extend 
(Sections \ref{facet-order-section} and \ref{uniqueness-section}) 
the applicability of the
methods from [BH] by using critical cell cancellation via gradient path reversal.
Most notably, we introduce the notion of a content-lex facet ordering, which has
also recently proven useful in work of [HHS].
The discrete Morse function on $\Delta (\Lambda) := 
\bigcup_{\lambda \in \Lambda}\Delta(\hat{0},\lambda)$ 
also yields a free resolution of $\field$ over $\field[\Lambda]$ whose 
multigraded Betti-numbers can be read off from the number of 
critical cells of given dimension and given multidegree, i.e. from 
the Morse numbers, and from the gradient paths governing incidence 
among critical cells.  Sections ~\ref{section-quad-applic} and 
~\ref{survive-deg-d} will describe the critical cells in the  
Morse function as follows: 
 
\begin{thm} 
The critical cells of the discrete Morse function on $\Delta (\Lambda )$  
in bijection with the words of a language accepted by a finite state  
automaton, i.e. the words of a regular language.  Thus, the 
generating function for Morse numbers is a rational function 
whose coefficients give upper bounds on all of the  
Betti numbers.   
\end{thm} 
 
In the quadratic Gr\"obner basis case, 
this generating function is exactly the Poincar\'e-Betti series,  
yielding yet another proof of its rationality in this case.  
There are known rational function bounds on the Poincar\'e-Betti
series (see Proposition 3.3.2 in \cite{Av}). But in contrast to these
bounds our rational function comes close to the actual Betti numbers, in the 
sense that it exhibits the vanishing of Betti numbers as in  
Theorem ~\ref{betti-vanish}. 
  
We conclude the paper in Section \ref{remarks-questions} with 
remarks and open questions. 

\section{Background} 

\subsection{Posets and Order Complexes} \label{posets}

Let $P$ be a poset with unique minimal element $\hat{0}$ and unique
maximal element $\hat{1}$. We denote by $\Delta(P)$ the simplicial complex
whose $i$-simplices are the chains $\hat{0} < p_0 < \cdots < p_i <
\hat{1}$ in $P$. The maximal chains -- with respect to inclusion -- 
in $P$ are sometimes called {\it saturated chains}.  
Notice that the  saturated chains in $P$ give rise to the facets (maximal faces) in  
$\Delta (P)$, and sometimes we will speak of saturated chains of $P$ and 
facets of $\Delta (P)$ interchangeably.   
For $x \leq y$ in $P$, let $[x,y]$ be the closed interval
$\{z~|~x \leq z \leq y\}$ and $(x,y) := [x,y] -\{x,y\}$ the open interval. 
We write $\Delta(x,y)$ for $\Delta([x,y])$.

For an arbitrary simplicial complex we denote 
by $H_i(\Delta;R)$ and $\widetilde{H}_i(\Delta;R)$ 
the non-reduced and reduced simplicial homology of $\Delta$ with coefficients 
in the ring $R$. 
The $i$-th Betti number $b_i$ of $\Delta$ (with respect to ${\bf Z}$) is 
the rank of the
free part of the $i$-th non-reduced homology group 
of $\Delta$ with coefficients in 
${\bf Z}$. 
In order to calculate the non-reduced homology $H_i(\Delta;R)$ or
the reduced homology $\widetilde{H}_i(\Delta;R)$ for the relevant simplicial
complexes $\Delta$ we will use two basic facts from algebraic topology.
First, if $\Delta$ is homotopy equivalent to a
topological space $X$ then the simplicial homology of $\Delta$ and the
cellular/singular homology of $X$ coincide. In our situation $X$ will always be a
CW-complex. The second important fact is, that if a CW-complex $X$ has 
$m_i$ cells of dimension $i$ then the 
$i$-th Betti number $b_i$ satisfies $b_i \leq m_i$.
A poset $P$ if called {\it homotopically Cohen-Macaulay} if the order complexes 
$\Delta(x,y)$ of all intervals $[x,y]$ in $P$  
are homotopy equivalent to a wedge of spheres of dimension $\dim(\Delta(x,y))$.
In particular this holds if for each $\Delta(x,y)$ either $\dim(\Delta(x,y)) \leq 0$ or
$\dim(\Delta(x,y)) \geq 1$ and $\Delta(x,y)$ is homotopy equivalent to a 
CW-complex with no cell in dimension $0 < i < \dim(\Delta(x,y))$ and a single cell in
dimension $0$.

Our main tool for the construction of a CW-complex $X$ homotopy equivalent to a given 
simplicial complex $\Delta$ is discrete Morse theory.

\subsection{Discrete Morse Theory: General Theory}
This section reviews discrete Morse theory results we will need 
from [Fo], [Ch], [Jo], [BH], [He2] and [BW], along with some  
other requisite background.  
Forman [Fo] defines a function $f$ which assigns real  
values to the cells in a regular 
CW-complex $X$ to be a {\it discrete Morse function} if  
for each cell $\sigma \in X^{(*)}$ the sets 
$$ \Big\{ \tau \subseteq \overline{\sigma}~\Big|~ \tau \in X^{(*)}, \dim(\tau) = \dim(\sigma)-1, f(\tau)  
\ge f(\sigma)\Big\}$$ and
$$ \Big\{ \overline{\tau}\supseteq \sigma 
~\Big|~ \tau \in X^{(*)}, \dim(\tau) = \dim(\sigma)+1, f(\tau) \le f(\sigma )\Big\}$$ 
each have cardinality at most one. Here $X^{(*)}$ denotes the collection of
open cells in $X$ and $\overline{\sigma}$ denotes the closure of 
$\sigma$ in $X$ for $\sigma \in X^{(*)}$. The condition implies that  
for each $\sigma$, at most one of two sets is non-empty. 
When both are empty, then $\sigma$ is called a {\it critical cell}. 
The main result on discrete Morse functions is the following:

\begin{thm}[\cite{Fo}] \label{main-discrete}
If $f$ is a discrete Morse function on the regular CW-complex $X$ then
$X$ is homotopy equivalent to a (not necessarily regular)  
CW-complex $X^M$, such that for any given $i$ the number of cells 
of dimension $i$ in $X^M$ equals the number of 
critical cells of dimensions $i$ of the Morse function of $f$.  
Moreover, incidences among cells in $X^M$ are 
governed by a collapsing procedure that leads from $X$ to  
$X^M$ while preserving homotopy type at each step.
\end{thm}

In [Ch], Chari 
reformulated discrete Morse functions for regular CW-complexes in terms of 
certain types of face poset matchings.  Recall that the {\it face poset} 
of a CW-complex $X$ is the partial order $F(X)$ on the cells in $X^{(*)}$ defined by  
$\tau \le \sigma$ whenever $\tau$ is contained in the closure 
$\overline{\sigma}$ of $\sigma$. If $X$ is the geometric realization of
an abstract simplicial complex $\Delta$ then this order is just the
inclusion relation between the simplices of $\Delta$.
The {\it Hasse diagram} 
of a poset is the graph whose vertices are the poset elements and whose 
edges are the covering relations $x\prec y$, i.e. pairs $x < y$ such that 
$x\le z \le y$ implies $z=x$ or $z=y$. 
 
\begin{defn}[\cite{Ch}] 
\rm{ 
A matching on the Hasse diagram of the face poset $F(X)$ of a regular CW-complex 
$X$ is called {\it acyclic} if the directed  
graph obtained by directing matching edges 
upward and all other poset edges downward has no directed  
cycles. 
} 
\end{defn}   
Notice that the  
non-critical cells of a discrete Morse function $f$ come in pairs that 
prevent each other from being critical. Hence, this pairing gives a matching on 
the face poset of the CW-complex. Furthermore, this matching is  
acyclic, because Chari's edge orientation will orient all edges  
in the direction in which $f$ weakly decreases.   
Conversely, 
many different (but in some sense equivalent) discrete Morse functions  
may be constructed from any face poset acyclic matching.  For instance, 
one may obtain $f$ by choosing a monotone function on any total order 
extension of the partial order given by the acyclic directed graph.  The  
face poset elements that are left unmatched by an acyclic matching 
are exactly the critical cells in any corresponding discrete Morse function. 
We will work exclusively in terms of acyclic matchings rather than  
discrete Morse functions, but at times it is helpful to have both points 
of view in mind. 
 
Denote by $m_i$ the number of critical cells of dimension $i$ in a discrete 
Morse function on a regular CW-complex $X$. As usual $b_i$ is the $i$-th
Betti number of $X$.
By virtue of Theorem \ref{main-discrete} $X$ is homotopy equivalent to 
a complex $X^M$ constructed from $m_i$ cells of dimension $i$.
The first of the following two results is an immediate corollary from
Theorem \ref{main-discrete}, the second was first proved in [Fo].
Both results exhibit a strong 
analogy with traditional Morse theory: 
\begin{eqnarray} 
m_j \ge b_j  \mbox{~for~}  0\le j\le \dim(X) \label{bound} \\ 
\displaystyle{\sum_{i=0}^{\dim (X)}} 
(-1)^i m_{\dim(X)-i} = \displaystyle{\sum_{i=0}^{\dim (X)}}  
(-1)^i b_{\dim (X )-i} = \chi (X ) 
\end{eqnarray}
 
\noindent The inequality (\ref{bound}) will be used in later sections in oder
to obtain bounds on Betti numbers. For our applications
also the following special situation which we already mentioned in
Section \ref{posets} will be of importance.
 
\begin{cor}\label{cm-ness} 
Let $X$ be a regular CW-complex and $f$ a discrete Morse functions on $X$ with
$m_i$ critical cells of dimension $i$.
If $m_i=0$ for $0 < i <j$ and $m_0 = 1$ then $X$ is $(j-1)$-connected. 
In particular, if for all order-complexes $\Delta(x,y)$ of intervals $[x,y]$ in a poset 
$P$ such that $\dim(\Delta(x,y)) > 0$ there is a Morse function $f_{xy}$ with no critical cell in
dimension $0 < i < \dim(\Delta(x,y))$ and a single critical cell in dimension $0$, 
then each interval is homotopy equivalent to a wedge of spheres of dimension $\dim((\Delta(x,y))$
and $P$ is homotopically Cohen-Macaulay. 
\end{cor} 
 
\begin{defn} 
\rm 
{ Let $X$ be a regular CW-complex and $f$ a discrete Morse function on $X$. 
A {\it gradient path} from a critical cell $\tau$ to another 
critical cell $\sigma$ of dimension $\dim(\tau)-1$ is a directed path upon 
which the Morse  function weakly decreases.} 
\end{defn} 

It is a simple consequence of the definition that a gradient path 
between $\tau$ and $\sigma$ will alternate between cells of
dimension $\dim(\tau)$ and $\dim(\tau)-1$.

It will turn out that for our purposes we need to find a
discrete Morse function on $\Delta(x,y)$ for intervals $x \leq y$ in 
an affine semigroup $\Lambda$ such that inequality (\ref{bound}) becomes an
equality. We will not always be able to achieve this goal and indeed it is
open whether this is even possible.
In our approach we construct a discrete Morse function on $\Delta(x,y)$ 
and then try to optimize the function in order to make (\ref{bound}) 
tight. For the latter we will employ the following observation.
 
\begin{obs}[\cite{Fo}] \label{cancelling} 
If there is a unique gradient path from $\tau$ to $\sigma$, then  
reversing the orientation of each edge in this path yields a  
new acyclic matching for which $\tau$ and $\sigma$ are no longer
critical.  
\end{obs} 

We refer to the procedure described in Observation \ref{cancelling}
as ``cancelling critical cells.'' 

The following lemma from [Jo] will  
be useful for combining several acyclic matchings to a single 
acyclic matching.
 
\begin{lem}[Cluster Lemma]\label{filter} 
Let $X$ be a regular $CW$-complex which decomposes into collections  
$X_p$ of cells indexed by the elements $p$ in a partial  
order $P$ with unique minimal element  
$\hat{0}$ as follows: 
\begin{enumerate} 
\item 
$X$ decomposes into the disjoint union $\cup_{p\in P}  
X_p$, that is,  
each cell belongs to exactly one $X_p$  
\item 
For each $p \in P$, $\cup_{q \le p} X_{p}$  
is a subcomplex 
of $X$ 
\end{enumerate} 
For each $p \in P$, let  
$M_{p}$ be an acyclic matching on the subposet $F(X) \cap X_p$ 
of $F(X)$ consisting of the cells in  $X_p$. 
Then $\cup_{p\in P} M_p$ is an acyclic matching  
on $F (X)$. 
\end{lem} 
 
\subsection{Discrete Morse Theory: The Case of Poset Order Complexes}
\label{poset-disc-review} 
 
Let $P$ be a poset with unique minimal element $\hat{0}$ and unique
maximal element $\hat{1}$. Assume that the edges of the Hasse diagram of $P$
are labeled by a labelling $\frl$ which takes values in a linearly
ordered set. Via this labelling we can assign to each saturated chain
$p_0 < p_1 < \cdots < p_{i-1}$ in $P$ of cardinality $i$ an $(i-1)$-tuple
$(\frl(p_0,p_1), \ldots, \frl(p_{i-1},p_i))$. Thus any linear extension
of the lexicographic ordering on the tuples induces a linear order on the
saturated chains in $P$.

In our situation where $P = [x,y]$ for $x \leq y$ in an affine
semigroup $\Lambda$ generated by $\alpha_1, \ldots, \alpha_n$ 
the labelling $\frl$ is given by sending a cover relation $\lambda \prec \mu$ 
to $\alpha_i = \mu - \lambda$ or if we consider $k[\Lambda]$ as a quotient of
$k[z_1, \ldots, z_n]$ we can equivalently label $\lambda \prec \mu$ with
$z_i$. Note, that in the latter case the product over the labels of a 
saturated chain is a monomial in $[z_1, \ldots, z_n]$.
In either case we can choose an arbitrary linear order on the sets
$\{\alpha_1, \ldots, \alpha_n\}$ or $\{z_1, \ldots, z_n\}$.
Since the two labellings are equivalent we will not distinguish between
$z_1$ and $\lambda_i$ in the rest of the paper.
 
It is shown in [BH] how to construct a discrete Morse function 
on $\Delta (P)$, respectively an acyclic matching on the face poset 
$F(\Delta(P))$ of  $\Delta (P)$ for a labeled poset $P$, from the 
lexicographic order on the saturated chains of $P$. If this lexicographic ordering 
happens to be a shelling order (see [BjW]) then it is possible to infer directly from the 
constructed Morse function that the poset is homotopically
Cohen-Macaulay -- just as it can be deduced from the lexicographic shelling itself.
We will refer to a Morse function resulting from a lexicographic ordering
on the saturated chains as a {\it lexicographic discrete Morse function}.  
We will show in  Section ~\ref{facet-order-section} that the construction 
of [BH] applies to a  larger class of orders on saturated chains.  This 
will allow us to construct a discrete Morse function on $\Delta(\Lambda)$ for
an affine semigroup $\Lambda$ generated by $\alpha_1, \ldots, \alpha_n$ from a  
facet order based on an arbitrary monomial term order on $k[z_1,\ldots, z_n]$.

Let $\frl$ be a labelling of the poset covering relations such 
that $\frl(u,v) \ne \frl(u,w)$ for $v\ne w$. Then the  ordering of 
the label sequences on the saturated chains by the lexicographic order 
already gives a total order on saturated chains, i.e. an ordering 
$F_1,\dots ,F_r$ on facets in $\Delta (P)$.  To describe the corresponding 
discrete Morse  function from [BH], we will need to speak of the ranks 
of elements of a saturated chain, whether or not a poset is graded. 
We will do this by assigning to an element of a saturated chain the rank of
the element within the chain and speak of the rank with respect to the chain.
Indeed, we do not require consistency of the notion of rank between 
different saturated chains. Also we will identify a set of chains in 
$P$, resp. faces of $\Delta(P)$, with the simplicial complex generated 
by the chains. 

Each maximal face in $F_j \cap (\cup_{i<j} F_i)$ has rank set 
(with respect to $F_j$) 
of the form $1,\dots ,i,j,\dots ,n$, i.e. it consists of all ranks in $F_j$ 
except for a single interval $i+1,\dots ,j-1$ of consecutive ranks that 
are omitted; this follows from the use of a lexicographic order on facets. 
Call each such list $i+1,\dots ,j-1$ of ranks a 
{\it minimal skipped interval} of $F_j$, and say the interval has  
{\it height $j-i-1$}.  Following [BiH] and [BH], call the collection of 
minimal skipped intervals for $F_j$ the {\it interval system} or {\it set of 
$I$-intervals} of $F_j$. 
 
For each facet $F_j$, [BH] constructs an acyclic matching on the 
set of faces in $F_j \setminus  
\cup_{i<j} F_i$ in terms of the interval system for  
$F_j$.
This is done in such a way that the  
union (over all $F_j$) of these matchings is acyclic, and 
each $F_j \setminus (\cup_{i<j} F_i)$ includes  
at most one critical cell.  We say $F_j$  
{\it contributes a critical cell} if $F_j \setminus (\cup_{i<j} F_i)$  
contains a critical cell.   
$F_j$ will contribute a critical cell if and only if the 
homotopy type changes with the attachment of $F_j$. 
 
\medskip  

\noindent {\bf Description of critical cells:}  
\begin{itemize} 
\item[(Case 1)] If the $I$-intervals of $F_j$ do not collectively have support covering 
all the ranks in $F_j$, then $F_j$ does not contribute a critical cell. 
\item[(Case 2)] If the $I$-intervals of $F_j$ cover all ranks and  
have disjoint support, then the critical 
cell consists of the lowest rank from each of the $I$-intervals. 
\item[(Case 3)] If there is some overlap in the minimal skipped intervals of $F_j$, then  
iterate the following procedure to obtain the critical cell, ordering  
$I$-intervals so that their minimal ranks are increasing: 
\begin{enumerate} 
\item Include the lowest rank from $I_1$ in the critical cell. 
\item Truncate all the remaining minimal. 
 skipped intervals by chopping off any ranks that they share with $I_1$.
\item Discard $I_1$ and any skipped intervals that are no longer minimal. 
\item Re-index the remaining truncated minimal skipped  
intervals to begin with a new $I_1$.
\item Repeat until there are no more minimal skipped intervals.
\end{enumerate} 
\end{itemize} 
  
The non-overlapping intervals obtained by the above 
truncation procedure are called the 
{\it $J$-intervals} of $F_j$. 
 
\begin{rk} 
$F_j$ contributes a critical cell if and only if its $I$-intervals  
cover all ranks of $F_j$.  In this case, the dimension of the critical 
cell is one less than the number of $J$-intervals. 
\end{rk} 
 
In order to cancel critical cells by reversing gradient paths, we 
will also need some  
information about the matching itself. 
 
\medskip  

\noindent {\bf Description of the acyclic matching on $F_j \setminus \cup_{i<j} F_i$:} 
\begin{itemize} 
\item If $F_j$ has no critical cells, then there is at least  
one cone point in $F_j \cap (\cup_{i<j} F_i)$. In this situation we match  
by including/excluding cone point of lowest rank.
\item If $F_j$ contributes a critical cell, then match any non-critical cell based on  
the lowest $I$-interval of $F_j$ where the cell 
differs from the critical cell.  Specifically, match by  
including/excluding the lowest element of the $I$-interval, since 
the cell must include at least one element of the $I$-interval other 
than this lowest possible element, in order to cover the $I$-interval 
but differ from the critical cell. 
\end{itemize} 
 
\begin{rk}\label{intuitive-idea-remark} 
If there is some $d$ such that every $I$-interval in a lexicographic 
discrete Morse function has height at most $d-1$, then the above  
construction immediately implies that each poset interval $(x,y)$ is 
at least $(-1 + \frac{\rm{rk}(y)-\rm{rk}(x)-1}{d-1})$-connected.  
In 
fact, it suffices for the average interval height to be at most $d-1$. 
\end{rk} 
 
In our setting, where $P = \Delta(x,y)$ for $x \leq y$ in an affine
semigroup $\Lambda$ Remark \ref{intuitive-idea-remark} will turn out
to be applicable for $d$ the degree of a Gr\"obner basis of the ideal
$I_\Lambda$ (i.e. the maximal total degree of a polynomial in the
Gr\"obner basis). We will use the fact that every $I$-interval results 
from a label sequence descent or from a syzygy leading term, by virtue of 
our use of a monomial term order to order the saturated chains. Furthermore, all  
leading terms will be divisible by Gr\"obner basis leading terms of 
degree at most $d$, which will allow us to show that the average interval 
height is at least $d-1$, yielding the desired lower bound on connectivity. 
In a sense, this gives a new combinatorial explanation for the  
connection between Gr\"obner basis degree and complexity, in which a  
shelling is the special case with $d=2$.  However, just as in the shelling 
in [PRS], only certain types of monomial term orders and Gr\"obner bases 
will immediately yield the degree $d$ analogue of a shelling.   

\begin{examp}\label{non-shell-examp} 
Consider $k[z_1,z_2,z_3,z_3]/(z_1z_4-z_2^2)$ with a term order
$\preceq$ such that $\init_\preceq(z_1z_4-z_2^2) = z_1z_4$. 
If $z_1z_4$ is a Gr\"obner basis leading term, but $z_1z_3$ and $z_3z_4$ 
are not Gr\"obner basis leading terms, 
Then $z_3(z_1z_4-z_2^2) = 0$ precludes the shelling from [PRS], since
$z_1z_3$ and $z_3z_4$ are not leading terms, 
\end{examp} 

In order to be able to allow completely general 
monomial term orders, i.e. to deal with situations such as  
in Example ~\ref{non-shell-examp}, we need to extend the 
tools from [BH] described in this section by performing critical 
cell cancellation.
This extension will be done in Section \ref{facet-order-section}. The the 
method of critical cell cancellation is explained in the 
next section. 
 
\subsection{Discrete Morse Theory: Optimizing  
Discrete Morse Functions}
 
This section reviews tools from [He2] for  
eliminating critical cells by cancelling pairs by 
a gradient path reversal.  Later we will construct 
a lexicographic discrete Morse function for monoid posets, and then 
use these Morse function optimization 
tools to eliminate all of the low-dimensional critical cells. 
 
Let $\Delta$ be a simplicial complex and $f$ a discrete Morse function on $\Delta$.
Define the {\it multi-graph face poset}, 
denoted $F(\Delta)^M$, for the complex $\Delta^M$ of critical cells as follows: 
\begin{enumerate} 
\item 
The vertices in $F(\Delta)^M$ are the cells in $\Delta^M$, or equivalently the  
critical cells in the discrete Morse function on $\Delta $. 
\item 
There is one edge between a pair of cells $\sigma ,\tau$ of  
consecutive dimension $\dim(\tau) = \dim(\sigma)+1$ for each gradient path from $\tau$ 
to $\sigma $. 
\end{enumerate} 
 
\begin{thm}\label{multi-graph} 
Any acyclic matching on $F(\Delta)^M$ specifies a collection 
of gradient paths in $F(\Delta )$ that may simultaneously be reversed to  
obtain a discrete Morse function $M'$ whose critical cells are the   
unmatched cells in the matching on $F(\Delta)^M$.   
\end{thm} 
 
To cancel cells, we will need to know that a 
gradient path from a critical cell $\tau $ to a critical cell 
$\sigma $ is the only gradient path from $\tau $ to $\sigma $.   
 
\begin{defn}\label{def-commute} 
\rm{ 
Let $u\prec v\prec w$ be covering relations in a poset
$P$ labeled by $\frl$. Assume we are given a linear order on the
saturated chains in $[u,w]$.
The labels $\frl(u,v)$ and $\frl(v,w)$ on 
covering relations $u\prec v\prec w$ 
are said to  {\it commute} if the  
least saturated chain in $[u,w]$ 
is labeled by $\frl(u,v) ,\frl(v,w)$ arranged in ascending order. 
} 
\end{defn} 
 
Let $\frl(u,v)$ denote the sequence of edge labels  
on the least saturated chain from $u$ to $v$.   
 
\begin{defn} \label{lci} 
\rm 
{ 
Let $P$ be a poset labelled by $\frl$ and assume that within each interval
the saturated chains are linearly ordered by an order depending only on
the label sequence.
\begin{itemize}
\item[(i)] The weakly increasing rearrangement of the label sequence of a
saturated chain is called the {\it content} of the chain.
\item[(ii)] The labelling $\frl$ is called {\it least-increasing} if every interval has a 
(weakly) increasing chain as its least saturated chain. 
\item[(iii)] The labelling $\frl$ is is called {\it least-content-increasing} if 
it is east-increasing and in addition the  
label sequence of the least chain equals or precedes the  
content of the label sequence of every other saturated chain in the interval. 
\end{itemize}
} 
\end{defn} 

Note, that since the definition assumes the linear order on the saturated chains only
depends on the label-sequences,  Definition \ref{lci}(iii) we can consider this order 
as an order on label sequences. If the linear order on the saturated chains is given by 
the lexicographic order then the condition least-increasing
is weaker than being an EL-labelling (see [BjWa]), in that intervals may have  
several increasing chains. 
 
\begin{rk} 
If two critical cells $\tau ,\sigma$ 
in a least-content-increasing labelling are contributed by saturated chains 
of equal content, then every downward step in any gradient path from 
$\tau $ to $\sigma $ must preserve content, and in fact must 
sort labels on the interval 
where the chain element was deleted. 
\end{rk} 
 
Combining results (see Theorem 6.6) in [He2] yields the following. 
 
\begin{thm}\label{red-exp} 
Let $P$ be a poset labelled by a least-content-increasing labelling $\frl$.
Let $M$ be the lexicographic discrete Morse induced by 
$\frl$ and let $\tau, \sigma$, $\dim(\tau) = 
\dim(\sigma)+1$, be critical cells resulting from saturated chains whose label 
sequences $\frl(\tau ), \frl (\sigma )$ have equal content.  
Suppose further that the permutation transforming  
$\frl(\tau )$ to $\frl(\sigma )$ is 321-avoiding.  
If there is a  gradient path $\gamma$  
from $\tau$ to $\sigma$ such that each downward step 
swaps a pair of consecutive labels by deleting an element $v$ from 
a chain which also includes elements covering and covered by $v$, then 
$\gamma$ is the unique gradient path from $\tau $ to $\sigma $. 
\end{thm} 
 
\begin{rk} 
Our upcoming discrete Morse function will use a content-lex facet 
order, as introduced in Section ~\ref{facet-order-section}.   
Theorem ~\ref{red-exp} also applies in that setting without requiring 
any modification. 
\end{rk} 
 
Theorem ~\ref{one-path} will  
generalize the above to deal  
with non-saturated chain segments, as needed for cancelling critical  
cells in our upcoming Morse function.   
 
\subsection{Discrete Morse Theory: Application to Cellular Resolutions}
\label{section-cellular} 
 
Let $M$ be a module over a commutative ring $R$. A {\it free resolution} 
of $M$ over $R$ is a complex of free $R$-modules $F_i$ and $R$-module   
homomorphisms $\partial_i$  
$${\mathcal F}: \cdots \mathop{\rightarrow}\limits^{\partial_i+1} F_i  
\mathop{\rightarrow}\limits^{\partial_i} 
F_{i-1} \mathop{\rightarrow}\limits^{\partial_{i-1}} \cdots  
\mathop{\rightarrow}\limits^{\partial_1} F_0.$$ 
which is exact in all degrees $\neq 0$ (i.e., $\Image (\partial_i )=  
\Ker (\partial_{i-1} )$ for $i \geq 2$) and $\Coker (\partial_1 )\cong M$. 
In our case $R = k[\Lambda]$ is a $k$-algebra and carries an  
additional multigraded structure,  
Recall that a $k$-algebra $R = k[z_1, \ldots, z_n] / I$ is called  
$\natnum^d$-{\it multigraded}   
if $R = \bigoplus_{\alpha \in \natnum^d} R_\alpha$  
as $k$-vector spaces and $R_\alpha  
R_\beta \subseteq R_{\alpha+\beta}$. If $d = 1$ and  
$R$ is generated in degree $1$  
over $k = R_0$ then $R$ is called {\it standard graded}. 
Analogously defined are $\natnum^d$-graded $R$-modules. In this  
situation we consider  
multigraded free resolutions ${\mathcal F}$. In addition to being  
a resolution 
one demands that the $F_i$ are free multigraded $R$-modules and that  
the $\partial_i$ 
are $\natnum^d$-homogeneous. A free multigraded module $F$ is a  
direct sum $\oplus_{\alpha \in \natnum^d} 
R(-\alpha)^{\beta_\alpha}$ of free $R$-modules $R(-\alpha)$ of rank one  
whose grading is defined 
by assigning $\alpha$ as the degree of the unit element $1$.  
 
A $\natnum^d$-graded free resolution  
$${\mathcal F}: \cdots \mathop{\rightarrow}\limits^{\partial_i+1}  
\bigoplus_{\alpha \in \natnum}  
R(-\alpha)^{\beta_{i,\alpha}} \mathop{\rightarrow}\limits^{\partial_i} 
\bigoplus_{\alpha \in \natnum} 
R(-\alpha)^{\beta_{i-1,\alpha}} \mathop{\rightarrow} 
\limits^{\partial_{i-1}} \cdots 
\mathop{\rightarrow}\limits^{\partial_1} \bigoplus_{\alpha \in \natnum} 
R(-\alpha)^{\beta_{0,\alpha}}.$$ 
is called cellular, if there is a CW-complex 
$X$ and a map $\gr : X^{(*)} \rightarrow \natnum^d$ from  
the set $X^{(*)}$ of its 
cells to $\natnum^d$ such that:  
\begin{itemize} 
\item There is a basis $e_c$ of $F_i = \bigoplus_{\alpha \in \natnum} 
R(-\alpha)^{\beta_{i,\alpha}}$ indexed by the $i$-cells of $X$  
in such a way that 
if $e_c$ belongs to $R(-\alpha)$ then $\gr (c) = \alpha$.   
\item For the $i$-cell $c$ of $X$ and its cellular differential  
$\delta_i (c) = \sum_{c'} [c : c'] c'$, we have  
$\partial_i(e_c) = \sum_{c'} [c : c'] \xx^{\gr (c') - \gr (c)} e_{c'}$. 
\end{itemize} 
In this situation we say that $X$ {\it supports} the  
resolution ${\mathcal F}$. 
 
Consider the face poset $P = F(X)$ of a CW-complex $X$ supporting a  
resolution.  [BW, Proposition 2.2] shows that an acyclic matching $A$ on $P$ 
leads to a chain homotopy between the original resolution and a smaller 
cellular resolution, given by the smaller CW complex $X^M$ of  
critical cells in a discrete Morse function given by $A$, if $A$  
matches only cells that have the same value under $\gr$. 
 
It is well known that for any $\natnum^d$-graded  
module $M$ there exists a {\it minimal} 
multigraded free resolution (i.e., a resolution that  
uses the least number of free 
modules in each degree). Now the results from \cite{BW}, as  
described above, 
allow one to construct smaller resolutions from a given  
resolution. It is also 
clear (see \cite{BW}) that this process will not always allow one to  
produce the minimal free resolution.  
 
Let us consider a `big' cellular resolution in the situation treated in this 
paper (i.e, $R = k[\Lambda]$ and $M = k$). 
It is well known that the simplicial complex $\Delta(\Lambda)$ of all finite 
chains $\lambda_0 < \cdots < \lambda_r$ in $\Lambda$ together  
with the grading 
$\gr(\lambda_0 < \cdots < \lambda_r) = \lambda_r$ gives a multigraded free  
cellular resolution of the maximal ideal  
$\umm = \Ker (k[\Lambda] \rightarrow k)$ -- 
the normalized Bar resolution. 
Since a free minimal resolution of $k$ over $k 
[\Lambda]$ starts with $k[\Lambda]$, minimizing the  
normalized Bar resolution  
is equivalent to minimizing a resolution of $k$. 
A well known criterion for a resolution to be minimal is that  
no unit elements of $R$ occur in the 
matrices representing the differentials. 
If in addition $R$ is standard graded   
and all matrix entries are either $0$ or elements of degree $1$, then the 
resolution is called {\it linear}. If $k$ has a linear resolution then 
$R$ is called {\it Koszul}.  
 
In our case, where $R = k[\Lambda]$, we know that $R$ is standard  
graded if and 
only if the generators of $\Lambda$ lie on an affine hyperplane.  
In this situation 
$R$ carries two gradings, the standard grading and a multigrading given by 
$\Lambda$. The Tor-groups $\Tor^R_i(k,k) \cong k^{\beta_i}$, where $\beta_i 
= \sum_{\alpha} \beta_{i,\alpha}$, also carry a multigraded structure  
$\Tor^R_i(k,k)_\alpha = k^{\beta_{i,\alpha}}$. Thus Koszulness  
can be read off from the $\Tor$-groups. Namely, $R$ is Koszul if and only if  
$\Tor^R_i(k,k)_\alpha = 0$ for $\xx^\alpha$ with standard grading not 
equal to $i$. 
 
A well known sufficient  
condition for a $k$-algebra $R =  
k[z_1,\ldots, z_n] /I$ which is 
standard graded to be Koszul is that $I$ has a quadratic Gr\"obner 
basis.  Recall that for a monomial order $\preceq$ on $k[z_1,\ldots, z_n]$ a  
{\it Gr\"obner basis} of $I$ is a generating 
set $\mathcal{G}$ of polynomials in $I$ 
such that the {\it initial ideal}  
$\init_\preceq (I) := \langle \init_\preceq (f) | f \in I\rangle$ is 
equal to the ideal generated by $\{ \init_\preceq (f) | f \in \mathcal{G}\}$. 
Recall that the 
{\it leading monomial}, denoted $\init_\preceq(f)$, for a polynomial $f$  
is the largest monomial with respect to $\preceq$ occurring in $f$; 
we write $\init(I)$ and $\init(f)$ if the monomial order is clear from the 
context.  Finally, a {\it monomial term order}  
on $k[z_1,\ldots, z_n]$ is a linear order $\preceq$ 
on the monomials in the ring such  
that (1) $1 \preceq m$ for all monomials $m$, and (2) $m \preceq m'$ implies 
$mn \preceq m'n$ for all monomials $n,m,m'$.   
 
\subsection{Basic Facts on Finite State Automata, Regular Languages and  
Rational Generating Functions}\label{regular-review} 
 
A central question in the theory of infinite resolutions is `Which conditions
on a module $M$ imply that Poincare'-Betti series of its minimal free resolution 
is rational ?' (see \cite{Av}). When $R = k[\Lambda]$ is a standard graded 
and multigraded $k$-algebra and $M = k$, the (graded and multigraded) {\it Poincare'-Betti series} is 
given by $\sum_{i,\alpha} \beta_{i,\alpha} t^i \xx^\alpha z^{\deg(\alpha)}$, 
where $\deg(\alpha)$ is the degree of $\xx^\alpha$ in the standard grading on 
$R$. 

We will not be able to give a new criterion for the rationality of the   
Poincare'-Betti series. But we will be able to give in Section \ref{survive-deg-d}
a rational series which bounds the Poincare'-Betti series from above,
i.e. all coefficients are greater or equal to the ones in the Poincare'-Betti series.
In order to prove rationality of our series we will resort to the theory of
regular languages. 
It is well known (see [BR]) that the generating series of a regular language is 
rational. A language $L$ over a finite {\it alphabet} $\Sigma$ is 
called {\it regular} if there is a finite state automaton which accepts 
exactly the words in $L$. The generating series of $L$ is given by 
$\sum_{w \in L} t^{|w|}$, where $|w|$ is the number of letters in $w$.  
See for instance [BR] for additional information. 

\section{Monomial term orders and discrete Morse functions 
resulting from (not-necessarily-lexicographic) facet 
orders}\label{facet-order-section}
 
This section will show how
the lexicographic discrete Morse function 
construction of [BH] generalizes easily
to a larger class of facet orders for poset order complexes; this
will include facet orders for monoid posets based on arbitrary
monomial term orders.  First observe that the 
[BH] construction applies without modification to any
facet order $F_1,\dots ,F_r$
which yields an interval system structure on each
$F_j \setminus (\cup_{i<j} F_i)$.  Equivalently,
the construction will work for facet orderings satisfying the crossing 
condition, as introduced in [He1] and defined below.  

{\bf Crossing condition.}  {\it Let $\leq$ be a linear order on the
saturated chains in a partially ordered set $P$ of rank $n$ and rank
function ${\rm rk}$. Let $F$ be a saturated chain,
$G \leq F$ and $\sigma = F \cap G$. Suppose 
that $[n] - \{ {\rm rk}(p)~|~p \in \sigma\}$
is not an interval of natural numbers. 
Then there is some facet $G' \leq F$ such that $F \cap G \subsetneq 
F \cap G'$.}
 
\medskip
The crossing condition implies that for a saturated chain $F$, 
maximal faces in $F \cap (\cup_{G < F} G)$ are supported on a set
of ranks whose complement is a single interval of consecutive ranks.

\begin{thm}\label{gen-lex}
If a facet ordering on an order complex satisfies the 
crossing condition, then the acyclic matching construction of [BH]
applies to this facet ordering.
\end{thm}

\proof
The effect of the crossing condition for a particular facet order 
$F_1,\dots ,F_k$ is to ensure that each maximal face of $F_j \cap (\cup_{i<j}
F_i)$ for $1 < j \le k$ skips a single interval of consecutive ranks.  
Thus, the faces in $F_j\setminus (\cup_{i<j} F_i)$ are the ones that ``hit''
each of these intervals, implying that the matching construction from
[BH] still applies.
\EOP

\begin{defn}
\rm
{
A facet order on a poset order complex which satisfies the crossing
condition is called {\it lex-like facet order}.  The discrete 
Morse function obtained by applying the construction from [BH] to such a 
facet order is called a {\it lex-like discrete Morse function}.
}
\end{defn}

Let $[x,v]$ be an interval in a labelled poset $P$. For a saturated chain in
$[x,y]$, the set of all saturated chains in $[x,y]$ having the same content is 
called the fibre of the content. Now again turn to intervals in affine
semigroups $\Lambda$. Let $\frl$ be the usual labelling of a covering relation 
$\lambda \prec \mu$ by the generator $\mu - \lambda$ of the semigroup.
We assume that we are given a monomial term order on $k[z_1, \ldots, z_n]$ and
again as usual identify the generators of $\Lambda$ with the variables $z_i$. 
This identification allows us to order the saturated chains by the given monomial 
term order. Recall, that we identify a label sequence with the monomial which is
the product over the labels, which in turn can be seen as the content of the
saturated chain with the given label sequence.
Lemma ~\ref{monoid-cross} considers linear 
orderings of saturated chains in intervals of $\Lambda$ obtained by combining the 
(commutative) monomial term ordering with a lexicographic order on each fibre; 
that is, we extend the order given by the monomial term order on the
content by the lexicographic ordering on chains that have the same content.
The lexicographic order uses the
monomial term order on degree $1$ monomials to order the labels.
 
\begin{lem}\label{monoid-cross}
Let $P = [x,y]$ be an interval in an affine semigroup $\Lambda$ and assume
that the saturated chains in $P$ are ordered by a monomial term order refined by
the lexicographic order. Then this ordering satisfies the crossing condition.  
In particular, we can construct a 
discrete Morse function just as in the case of a lexicographic order.
\end{lem}
 
\proof
Let $\equiv$ be the equivalence relation on the set of saturated chains
such that $m \equiv n$ if and only if $m$ and $n$ are two saturated chains
in the same closed interval in $P$.  We abuse notation and say
$m=n$ if the two chains are labeled by the same commutative monomial,
i.e. the labels on one saturated chain are a permutation of the labels on
the other.
Suppose that for two saturated chains $m_1m_2$ and $n_1n_2$ we 
have $n_1 \neq m_1$, $n_2 \neq m_2$ but $n_1 \equiv m_1$ and
$n_2 \equiv m_2$. Assume further that $n_1m_1 <_{monom} n_2m_2$ in the
given monomial order $<_{monom}$. 
This is one situation where a saturated chain has an overlap face with 
earlier
saturated chains such that the complement of the ranks in the overlap face
is disconnected.  
 
We check that either $m_1n_2 <_{monom} n_1n_2$ or
$n_1m_2 <_{monom} 
n_1n_2$, as follows.  Suppose $m_1n_2 >_{monom} n_1n_2$, which 
implies $m_1 >_{monom} n_1$.  Suppose 
$n_1m_2 >_{monom} n_1n_2$ also holds, implying
$m_2 >_{monom} n_2$.  Combining these inequalities yields
$$m_1m_2 >_{monom} n_1m_2 >_{monom} n_1n_2,$$ a contradiction.  Hence, 
at least one of the monomials $m_1n_2$ or $n_1m_2$ precedes
$n_1n_2$ in our term order.  This ensures that the maximal face shared
by $m_1m_2$ and $n_1n_2$ is not a maximal face in the simplicial complex of
faces shared by $n_1n_2$ and facets that precede it in lexicographic order,
just as needed.  

Now suppose there is a saturated chain $F_j$ not dealt with above that
has an overlap face 
with an earlier saturated chain $F_i$ such that the 
complement of the ranks in 
the overlap face is disconnected.  Then $F_i,F_j$ are labeled 
$m_1m_2m_3, \pi(m_1)n_2\sigma(m_3)$, respectively, where $\pi,\sigma $
are permutations on the labels in $m_1,m_3$ and $m_2\equiv n_2$ but
$m_2\ne n_2$.  Consider $F_i'$ which is labeled 
$\pi(m_1)m_2\sigma(m_3)$.  Since $F_i\cap F_j \subsetneq F_i'\cap F_j$,
we are done.
\EOP

\begin{rk}
\rm{
Examples suggest that 
the following procedure may be convenient for posets with no
particularly well-behaved global labelling.
\begin{enumerate}
\item
Label edges in a poset Hasse diagram in a natural way (or more generally,
give a chain-labelling on saturated chains).
\item
Partition the set of saturated chains into groups called {\it content
classes} according to the 
content of their label sequences.
\item
Put an ordering on these content classes.
\item
Within each content class, order saturated chains lexicographically.
\item
Prove that the resulting facet order satisfies the crossing condition.
\item
Cancel pairs of critical cells which have the same content, i.e. pairs
in the same content class.
\end{enumerate}
This approach
can make critical cell cancellation manageable 
for posets with no particularly nice global 
labelling, because gradient paths that begin and end in the same content
class must never leave that content class.  This will be essential to our 
analysis of monoid posets and to arguments in [HHS].
}
\end{rk}

\begin{defn}
\rm
{
A {\it content-lex facet order} is a lex-like facet order such that:
\begin{itemize}
\item The ordering is constructed from a labelling by refining a linear order on 
fibres by a lexicographic order.  
\item The least saturated chain in each interval has weakly increasing labels.
\end{itemize}
}
\end{defn}

One way content-lex facet orders arise is when
each content class individually has an EL-labelling.  This
is the situation for our upcoming facet order on monoid posets as well as
a $GL_n(q)$-analogue of the partition lattice examined in [HHS].

\begin{rk}
Content-lex facet orders behave as least-content-increasing
labellings for purpose of applying Theorem ~\ref{red-exp}
to verify gradient path uniqueness.  
\end{rk}

\section{Uniqueness of 321-avoiding gradient paths in 
discrete Morse functions from content-lex 
facet orders}\label{uniqueness-section}

Now we generalize Theorem ~\ref{red-exp} to allow non-saturated chain
segments in gradient paths between critical cells in the same content
class in a content-lex facet order, under certain additional 
assumptions.  Remark ~\ref{non-unique-remark} suggests 
that the 321-avoiding assumption is
probably necessary for any general result about gradient path uniqueness
for lex-like discrete Morse functions.

Upcoming sections will construct a non-optimal 
discrete Morse function for monoid posets, then use the following 
theorem to improve it.

\begin{defn}
\rm
{
A {\it delinquent chain} in a least-increasing labelling (or in a 
content-lex facet order) is an
increasing chain that is not lexicographically smallest on an
interval.  A content-lex facet order is {\it consistently delinquent}
if the existence of a delinquent chain labeled 
$a_1\cdots a_k$ implies that any chain segment 
labeled $b_1\cdots b_{k'}$ for $b_1=a_1,b_{k'}=a_k$ and 
$\{ a_1,\dots ,a_k\} \subseteq \{ b_1,\dots ,b_{k'} \} $ is
also delinquent.
}
\end{defn}

Denote by $e\alpha $ the earliest facet containing a cell $\alpha $,
and let $\frl(\alpha )$ be the label sequence on $e\alpha $.
We say that a chain element {\it covers} a delinquent
chain if it is in the interior of the delinquent chain, preventing the
poset chain from belonging to an earlier content class.  

\begin{thm}\label{one-path}
Let $\tau,\sigma$, $\dim(\tau) = \dim(\sigma)+1$, be critical cells
in the same content class in a content-lex facet order which 
is consistently delinquent.  Suppose $\frl(\sigma )$ differs from 
$\frl (\tau )$ by a 
321-avoiding permutation $\pi $ which
either shifts a single group of one or more 
consecutive ascending labels upward, or shifts a single label 
downward.  Then there is at most one gradient path 
from $\tau $ to $\sigma $.
\end{thm}

\proof
Since $\tau $ is critical, every 
pair of consecutive labels in $\frl(\tau )$ is 
either a descent or part of a minimal delinquent chain.  
Gradient paths from $\tau $ to $\sigma $ can never uncover any delinquent
chains, since that would cause the gradient path to pass
to an earlier content class, from which it could never reach
$\sigma $.  Since the facet order is also least-increasing, each downward 
step must either preserve the label sequence or 
eliminate a descent by deleting an element 
$v_r$ from a chain $\tau_i $ of the form $v_1 < \cdots < 
v_{r-1} < v_r < v_{r+1} < \cdots < v_s$, causing the labels on the segments 
of $e\tau_i $ from $v_{r-1}$ to $v_r$ and from $v_r$ to $v_{r+1}$ to be
sorted into a single ascending list.  The least-increasing property
Thus, any inversions
present in $\frl(\sigma )$ must be preserved throughout the gradient
path, since they can never be re-introduced.
With these observations in hand, we will describe the only possible 
gradient path from $\tau $ to $\sigma $, showing at each stage there
is only one choice for how to proceed.

Suppose $\frl(\sigma )$ is obtained from $\frl(\tau )$ by 
shifting a collection of consecutive, ascending labels upward.
Let $\mu $ be the label to be shifted upward to the 
highest destination in $\frl(\sigma )$, and let $b$ be the label
immediately above $\mu $ in $\frl(\tau )$.
By virtue of the Morse function construction of [BH], $\tau $ includes 
exactly the ranks where 
$e\alpha $ has descents as well as exactly one rank that covers each of
the minimal delinquent chains, namely the lowest ranks in the 
$J$-intervals.  The first gradient path downward step must eliminate 
a descent since it is not allowed to uncover a delinquent chain.  
The only choice that will not eliminate an inversion that is present
in $\frl(\sigma )$ is to eliminate a descent between 
$\mu $ and the label above it, i.e. deleting an element  
$v_r$ from a chain $v_1 < \cdots < v_{r-1} < v_r < v_{r+1} < \cdots 
< v_s$.  This yields a chain $\sigma_1 $ with $e\sigma_1 = F_j$ 
having a cone point in $F_j \cap (\cup_{i<j} F_i)$ between $v_{r-1}$
and $v_{r+1}$ in $F_j$, but 
no lower cone points.  Thus, the subsequent
upward step must insert some $v_{r'}$ above $v_{r-1}$ and below the label
$\mu $, with the labels between $v_{r-2}$ and $v_{r'}$ now 
comprising either a descent or minimal delinquent chain.  In the latter
case, the delinquent chain just below $v_{r'}$ must include the label
$b$ as its highest label.
Labels now below 
$b$ can no longer shift upward, since any such labels to be shifted 
upward in $\frl(\sigma )$ must be smaller than $b$.
If there were instead a descent at $v_{r-1}$, only one label is 
allowed between $v_{r-1}$ and $v_{r'}$, so either it is the label $b$ 
defined above, which does not shift upward, or it is  
a label smaller than $b$, which is now prevented from moving upward
by virtue of being smaller than $b$.

Continuing in this fashion, there is only one viable
downward step at any given stage until $\mu $ reaches its destination, 
namely the chain deletion which shifts $\mu $ upward,
since no $J$-interval is ever covered by more than one chain element
while $\mu $ is shifting upward.
By the argument above, 
labels can only shift upward while they are shifting as part of a block 
of consecutive labels which includes $\mu $, so all 
label shifting is complete once $\mu $ has reached its
destination, and there is a unique way for this to happen.  All that
remains is to consider additional gradient path steps which preserve label
sequence.

Suppose the label above $\mu $ in $\frl(\sigma )$
is larger than $\mu $.  Then $\mu $ must be within a delinquent chain,
since $\sigma$ is critical. 
After inserting a cone point
below $\mu $, the gradient path must
take a downward step deleting a chain element above $\mu $, since
the $J$-interval for the delinquent chain which includes $\mu $ will
now be covered by at least two chain elements.
At this point, no upward step is possible, so the gradient path must have 
reached $\sigma $.  
Finally, we show that $\mu $ cannot form a descent with the 
label above it in $\frl(e\sigma )$.  Otherwise, $\mu $
would be 
the highest label in a minimal
delinquent chain labeled $\nu_1,\dots ,\nu_k ,\nu,\mu $, but $\mu $ 
could not be shifted upward to this position from below without passing
to an earlier content class.  

The case where a single label $\mu $ shifts downward to obtain 
$\frl(\sigma )$ from $\frl(\tau )$ is quite similar, so is
essentially left to the reader.
The first gradient path step must again eliminate a descent so as to 
preserve content class, and the only allowable choice is a descent 
between $\mu $ and the label below it.  Similar reasoning to above will
show that the only possible gradient path will progressively shift 
$\mu $ downward to its destination, concluding once $\mu $ reaches the
interior of a delinquent chain.
\EOP

\section{The Cohen-Macaulay property for monoid posets with 
quadratic Gr\"obner bases}\label{cm-section}

Throughout this section, we assume the toric ideal
$I_{\Lambda } = \Ker(\phi )$ has a quadratic Gr\"obner basis $B$.  
However, Section ~\ref{non-optimal-section} will apply to higher 
degree Gr\"obner bases with essentially no modification needed,
and large parts of Sections ~\ref{syzygy-section}--~\ref{cm-subsection} 
will also generalize easily to higher degree Gr\"obner bases.

\begin{defn}
Denote
by $\init(I_{\Lambda })$ the initial ideal of $I_{\Lambda }$ with respect
to the term order giving rise to $B$, i.e. the ideal generated by 
leading terms of elements of $B$.  
\end{defn}

This section will
show that each interval $(\hat{0},m)$ in the resulting
monomial divisibility poset
has the homotopy type of a wedge of spheres of top dimension.  Our 
approach will be to construct a lex-like discrete Morse function based
on a content-lex facet order in
Theorem ~\ref{morse-init}, 
and then to cancel all but some of the top-dimensional critical cells in
Theorem ~\ref{morse-final}.  

We will employ the simple fact that any leading term of a polynomial in 
a toric ideal is divisible by a Gr\"obner basis leading term. In particular, 
leading terms of degree greater than the degree of the Gr\"obner basis
will include variables that in some sense are non-essential.
We will use these non-essential variables to cancel critical 
cells.  In the case of quadratic Gr\"obner bases, this approach will 
allow us to cancel all
critical cells that are not saturated chains. 
Later, we will 
use a similar (but somewhat more intricate) analysis for Gr\"obner
bases of degree $d$.

\subsection{A non-optimal Morse function}\label{non-optimal-section}

The first step will be to give a content-lex facet ordering. 

\begin{thm}\label{morse-init}
The monoid poset interval $(\hat{0},m)$ has a 
lex-like discrete Morse function resulting from a content-lex
facet order.  Its minimal skipped intervals
are the saturated chain segments with label sequences
of the following two types:
\begin{enumerate}
\item
descents
\item
sequences of weakly increasing labels
$\lambda_i,\dots ,\lambda_k $ such that 
$\lambda_i\lambda_k\in \init(I_{\Lambda })$, but 
$\lambda_{i'}\lambda_j\not\in \init(I_{\Lambda })$ for every other pair
$i\le i'<j\le k$. 
\end{enumerate}
\end{thm}

\proof
The finite saturated chains on intervals in a monoid poset 
correspond naturally to pairs $(m_i,\pi )$ where $m_i$ is a monomial in 
$k[z_1,\dots ,z_n]$, and $\pi $ is an ordering on the content of $m_i$.
This is equivalent to labelling saturated chains by 
non-commutative monomials, in $k\langle y_1,\dots ,y_n\rangle $, the 
viewpoint taken in [PRS].  Following [PRS], we order
saturated chains in an interval $(\hat{0},m)$ by using the
monomial term order which led to the Gr\"obner basis $B$ to order the
factorizations $m_i\in \phi^{-1}(m)$, and then 
lexicographically ordering label sequences of any fixed content, with 
our label order given by the monomial term
order applied to monomials of degree one.
Lemma ~\ref{monoid-cross} confirms the crossing condition for this 
facet order, implying it gives rise to a lex-like
discrete Morse function from a content-lex facet order.
Next, we characterize its minimal skipped intervals.  

Notice that a descent on the saturated chain segment $u\prec v\prec w$
implies a lexicographically smaller ascend $u\prec v'\prec w$, obtained
by reversing the order in which semi-group generators are multiplied, so 
descents always give minimal skipped intervals.  
On the other hand, any label sequence 
$a_1,\nu_1,\nu_2,\dots ,\nu_r ,a_2$ as in (2) will 
give rise to a minimal skipped interval because the 
Gr\"obner basis element $a_1a_2 - b_1b_2$
with leading term $a_1a_2$ implies the existence of an earlier saturated 
chain on the interval labeled by the increasing rearrangement of
the monomial $b_1n_1\cdots n_jb_2$; its minimality follows from the lack
of descents and of Gr\"obner basis leading terms not requiring both 
$a_1$ and $a_2$.  

\begin{defn}
\rm
{
The second type of minimal skipped
interval in the statement of the theorem
is called a {\it syzygy interval}.
}
\end{defn}

To see there are no other minimal skipped intervals, note that a label 
sequence $a_1,\dots ,a_r$ on any other minimal skipped interval must be
ascending to avoid descents which would preclude its minimality; to have
an earlier saturated chain on the interval 
$a_1,\dots ,a_r$ must be a leading term, hence divisible 
by a Gr\"obner basis leading term $m$.  But 
minimality ensures $a_1,a_r$ divide $m$, and the 
fact that the Gr\"obner basis is quadratic implies $m=a_1a_r$.
\EOP

\begin{cor}
Theorem ~\ref{one-path} may be applied to the above Morse function
to cancel pairs of critical cells that belong to the same fibre.
\end{cor}

\begin{defn}
\rm
{
A minimal skipped interval is {\it non-trivial} if it
has height greater than one.  Notice that only syzygy intervals may 
be non-trivial.
}
\end{defn}

\begin{examp}\label{lex_morse_examp}
\rm{
Consider the interval $(1,x_1^2x_2^2x_3x_4)$, or equivalently, 
$(1,z_1z_2z_3z_4) $, in the ring $k[x_1x_2,x_1^2,x_3,x_4,x_2^2] \cong 
k[z_0,z_1,z_2,z_3,z_4]/(z_1z_4-z_0^2)$, with $\init(I_{\Lambda }) = 
(z_1z_4)$.  Saturated chains are labeled by indices of the generators
$z_0,\dots ,z_4$.  
\begin{figure}[h]
\begin{picture}(250,100)(-45,0)
\psfig{figure=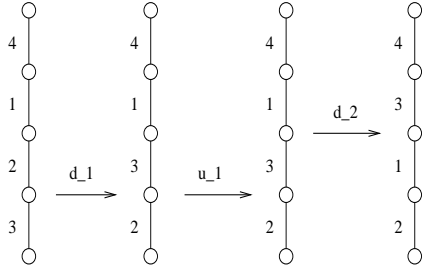,height=3.5cm,width=5.5cm}
\end{picture}
\caption{A gradient path shifting the label 3 into a syzygy interval}
\label{des_shift}
\end{figure}
Figure ~\ref{des_shift} shows four saturated chains on this interval, the
leftmost and rightmost of which will contribute critical cells.
Notice that $1,4$ labels a syzygy interval, while there are
descents at ranks 1 and 2 in the leftmost saturated chain, so it
has a critical cell $\tau $ comprised of 
ranks 1,2,3.  On the 
other hand, the label sequence $1,3,4$ in the rightmost saturated chain also
labels a syzygy interval, and this chain has a
descent at rank 1, so it has a critical cell $\sigma $ comprised 
of ranks 1,2.  
}
\end{examp}

\medskip
The remainder of Section ~\ref{cm-section} is 
devoted to cancelling pairs of critical cells 
by gradient path reversal, so as to eliminate all critical cells not
given by saturated chains. 
This will require an acyclic matching on $F(\Delta)^M$ consisting of pairs 
of critical cells to be cancelled.

\begin{defn}
\rm
{
We call the critical cells that remain after all this cancellation
the {\it surviving critical cells}.  
}
\end{defn}

\subsection{Syzygy intervals and their 
non-essential sets}\label{syzygy-section}

For now we assume all monomials on our monoid poset interval
are square-free.
The general case is dealt with in Theorem ~\ref{morse-final}.

\begin{rk}
For convenience,
we will refer interchangeably
to a critical cell and the saturated chain which contributes it.
\end{rk}

Figure ~\ref{des_shift} gives an example of a gradient path from a 
critical 2-cell to a critical 1-cell
in the order complex resulting from the semigroup
ring $k[z_0,z_1,z_2,z_3,z_4]/(z_1z_4 - z_0^2)$.  This gradient path
shifts the label $z_3$ to the interior of a syzygy interval, using the
fact that $z_3$ is not essential to $z_3(z_1z_4 - z_0^2) = 0$.
This is the only gradient path between these two critical cells.
Our goal will be to systematically cancel many such pairs
of critical cells simultaneously.

\begin{defn}
\rm
{
Denote by $I(a_1,a_2)$ the syzygy interval with ascending labels
$a_1,\lambda_1,\dots ,\lambda_k ,a_2$ in a saturated chain,
and refer to the Gr\"obner basis leading term $a_1a_2$ with $a_1\le a_2$
as an {\it increasing leading term}, or {\it ILT} for short.  
}
\end{defn}

We will soon use
ILTs to
collect critical cells into Boolean algebras within $F(\Delta)^M$.  

\begin{examp}
\rm{
In the affine semi-group ring $k[z_1,\dots , z_6]/(z_2z_6-z_1^2)$, 
consider the saturated chain 
$F_j$ that is labeled $z_4z_3z_2z_5z_6$.  $F_j$ contributes the 
critical cell $\sigma = z_4 < z_3z_4 < z_2z_3z_4$.
By Theorem ~\ref{red-exp}, there is a unique gradient path
from the critical cell
$\tau = z_5 < z_4z_5 < z_3z_4z_5 < z_2z_3z_4z_5$ to $\sigma $, 
given by the reduced expression $s_1\circ s_2\circ s_3$.  
More generally,  
each $T \subseteq S= \{ z_3,z_4,z_5\} $ 
gives rise to a critical cell $Crit(T)$ contributed by a facet $F_T$, 
as follows.  $F_T$ has label sequence $z_{(S\setminus
T)^{rev}}z_2 z_T z_6  $, where 
$z_T$ is the list of members of $T$ in increasing
order, and $z_{(S\setminus T)^{rev} }$ is the list of members of $S\setminus
T$ listed in decreasing order.  $S = \{ z_3,z_4,z_5\} $ is 
the non-essential set
of the interval.
Theorem ~\ref{red-exp} will show that the set of critical 
cells $\{ Crit(T)|T\subseteq S\} $ sits 
inside the multi-graph face poset $F(\Delta)^M$ as a Boolean algebra, depicted in
Figure ~\ref{fig-boolcrit}.  
\begin{center}
\begin{figure}[h]
\begin{picture}(250,360)(-50,10)
\psfig{figure=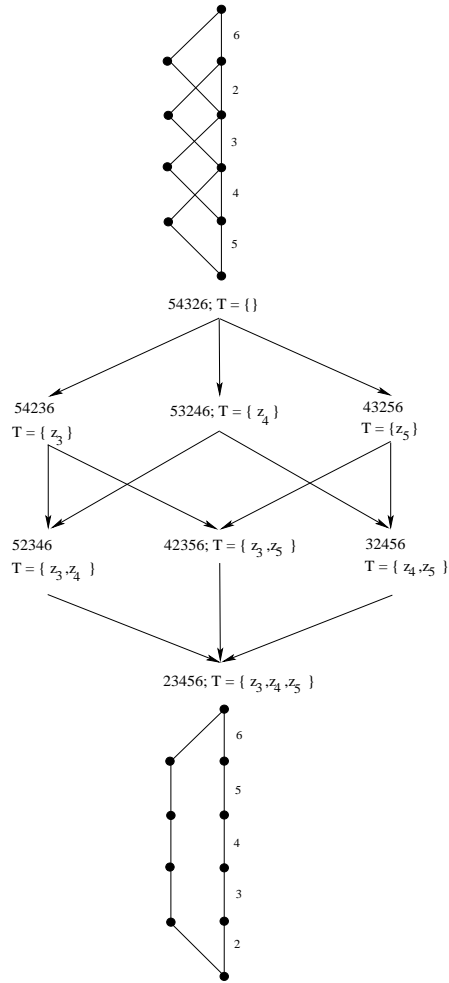,height=13cm,width=6cm}
\end{picture}
\caption{A Boolean algebra in $F(\Delta)^M$ indexed 
by subsets of $\{ z_3,z_4,z_5\}$ }
\label{fig-boolcrit}
\end{figure}
\end{center}
This Boolean algebra has covering relations $Crit(T\cup \{ z_i\} )\prec 
Crit(T)$ for each $T\subseteq S$ and each $ z_i\in S\setminus T$.  
}
\end{examp}

\begin{rk}
A gradient path cannot swap non-commuting labels (in the sense of 
Definition ~\ref{def-commute}) without passing to an 
earlier fibre, so cells to be cancelled will agree up to allowable
label commutation.
\end{rk}

\begin{defn}
\rm{
A label $\lambda $ in a saturated chain $M$
is {\it upward-shiftable} into a syzygy interval
$I(a_1,a_2)$ if it satisfies all the following conditions:
\begin{enumerate}
\item
$a_1 <_{monom} \lambda <_{monom} a_2$
\item
$\lambda $ appears below $I(a_1,a_2)$
\item
all labels between $\lambda $ and $I(a_1,a_2)$ are smaller than $\lambda $
and commute with $\lambda $ 
\item
all labels within $I(a_1,a_2)$ commute with $\lambda $
\item
$\lambda $ is not 
the top of some $I(\mu ,\lambda )$ with either non-empty interior or
such that the label $\nu $ immediately above $\lambda $ would 
neither form a descent with $\mu $ nor be part of an ILT together 
with $\mu $
\end{enumerate}
Likewise, $\lambda $ is 
{\it downward-shiftable} from $I(a_1,a_2)$ to just above $\lambda '$ if
$\lambda $ commutes with all labels separating it from $\lambda '$ and 
is larger than all such labels.
}
\end{defn}

\begin{rk}
Each label has at most one syzygy interval into which it
is upward-shiftable, because the lowest such interval
will separate it from all higher
ones.
\end{rk}

\begin{defn}
\rm
{
If $\lambda $ appears within $I(a_1,a_2)$, then 
the {\it topologically decreasing position} for $\lambda $ below 
$I(a_1,a_2)$ is the 
highest position below $I(a_1,a_2)$ to which $\lambda $ is
downward-shiftable so as to obtain the label sequence for
a critical cell with $\lambda $ not in the interior of any ILT,
if such a position exists.    
}
\end{defn}

In its topologically decreasing 
position below an ILT, $\lambda $ must form descents or ILTs
with the labels above and below it.
Lemma ~\ref{bool} will construct gradient 
paths that shift labels from their topologically decreasing positions
below ILTs to the interior of ILTs.
In some circumstances, we will also 
speak of the topologically decreasing position of $\lambda $ 
above $I(a_1,a_2)$, by which we mean the lowest position above 
$I(a_1,a_2)$ to which $\lambda $ is upward-shiftable to yield the 
label sequence for a critical cell.

\begin{defn} 
\rm
{
The {\it non-essential set}
of a syzygy interval $I(a_1,a_2)$ that appears in the label sequence
$\frl(\sigma )$
for a critical cell $\sigma $
will be a collection of labels that
appear in  
$\frl(\sigma )$ either
within $I(a_1,a_2)$ or 
in topologically decreasing positions below $I(a_1,a_2)$.
This set of labels (to be defined precisely in the remainder of 
Section ~\ref{syzygy-section} and Section ~\ref{match-subsection}) 
is denoted $S(a_1,a_2)$.
}
\end{defn}

To try to convey the intuition for $S(a_1,a_2)$, we now give an 
oversimplified definition.
Section ~\ref{match-subsection} will modify this into a 
much more technical definition that accomplishes exactly what is 
needed.
Initially, let us include in $S(a_1,a_2)$ those labels that appear within
$I(a_1,a_2)$ that are downward-shiftable to topologically decreasing
positions below $I(a_1,a_2)$.
Denote these labels 
by $n_1,\dots ,n_j$.
Also include in $S(a_1,a_2)$ those labels $m_1,\dots , m_r$ that are
upward-shiftable into $I(a_1,a_2)$, 
chosen in order from highest to lowest topologically 
decreasing position below $I(a_1,a_2)$.

We call 
critical cells which are not maximal faces in the order complex
$\Delta (\hat{0},m)$ {\it unsaturated}.
In the case of a quadratic Gr\"obner basis, we will match and 
cancel all unsaturated critical cells, using the fact that each
must have one or more syzygy intervals with non-empty interior.

\subsection{The matching on critical cells}\label{match-subsection}

In this section, we precisely define
non-essential sets and show that the resulting 
matching on critical cells is well-defined.
We also show that every critical cell which has at least one
syzygy interval with non-empty interior is indeed matched and 
cancelled, by showing it has at least one syzygy interval
with non-empty non-essential set.
The fact that pairs of critical cells to be matched 
do indeed comprise covering
relations in $F(\Delta)^M$ will be verified in Section ~\ref{cm-subsection}.

\begin{defn}
\rm
{
A label $\lambda $ is {\it preferable} to a label $\mu $ within a 
label sequence if either $\lambda $ is in the non-essential set of a 
higher syzygy interval than $\mu $ is in, or $\lambda ,\mu \in 
S(a_1,a_2)$ with the topologically decreasing position for $\lambda $
higher than for $\mu $.  
}
\end{defn}

Using our oversimplified definition of non-essential set from the 
previous section, 
let $I(a_1,a_2)$ be the highest syzygy
interval in a saturated chain $C$ such that $S(a_1,a_2)\ne \emptyset $,
and let $\lambda $ be the label in $S(a_1,a_2)$ with highest topologically
decreasing position below $I(a_1,a_2)$.
The theorem below will sometimes include in a non-essential set
$S(a_1,a_2)$ a single label that shifts 
downward into $I(a_1,a_2)$ from above. 
When this happens, denote this label $m_0$, and eliminate from
$S(a_1,a_2)$ any labels below $I(a_1,a_2)$ that do
not commute with $m_0$.

\begin{thm}\label{well-def-bool}
Every critical cell with a syzygy interval with an internal label
$\mu $ is matched.  Specifically, if
$\mu $ is excluded from the non-essential set of the interval, 
then there must be another label 
that allows the cell to be matched and cancelled.  Moreover, the
matching choices are made consistently.
\end{thm}

\proof
We will typically match by shifting $\lambda $ from inside
$I(a_1,a_2)$ to its topologically decreasing position below
$I(a_1,a_2)$, or vice versa.
However, special care is needed in four circumstances described below.
When $\lambda $ is excluded from a non-essential set, then the
matching instead shifts the label with highest preferability among those
belonging to some non-essential set.

In each of these circumstances, we will show that either $\lambda $ may be
included in $S(a_1,a_2)$ or that there is an alternative label to
$\lambda $ allowing the cell to be matched.
Moreover, when $\lambda $
is excluded from a non-essential set, it will be excluded for all
critical cells in the Boolean algebra within which our critical cell 
is matched.
\begin{enumerate}
\item
$\lambda \in I(a_1,a_2)$ shifts downward to a topologically decreasing
position which is higher than some label $\mu $ with which $\lambda $
does not commute; however, $\mu $ may
shift upward to the interior of a syzygy interval $I(b_1,b_2)$ to obtain
another critical cell when $\lambda $ appears in $I(a_1,a_2)$ but not 
when $\lambda $ has shifted downward to its topologically decreasing
position.
\item
$\lambda \in I(a_1,a_2)$ shifts downward to just above
a label $\mu $ with which $\lambda $ does not commute to form an
ILT $I(\mu ,\lambda )$ which then has 
some label $\nu \in S(\mu ,\lambda )$.
\item
$\lambda \in I(a_1,a_2)$ cannot shift downward to a topologically
decreasing position without first encountering a 
label $\mu $ with which it does
not commute.
\item
$\lambda \in I(a_1,a_2)$, but shifting $\lambda $ downward causes
$a_2$ to be in the non-essential set of a higher ILT, or more generally
the shifting of all labels within $I(a_1,a_2)$ downward to topologically
decreasing positions, cumulatively causes $a_2$ to belong to a higher
non-essential set.
\end{enumerate}

In the first case above, exclude
$\lambda $ from $S(a_1,a_2)$.  When $\lambda $ appears within
$I(a_1,a_2)$, then $\mu \in S(b_1,b_2)$, ensuring
the cell may still be matched.
Matching by 
shifting $\mu $ will clearly give a cell which would also exclude 
$\lambda $ from $S(a_1,a_2)$.  
Shifting a label $\mu '$ that is preferable to $\mu $ also
gives a cell that excludes $\lambda $, either by virtue of $\mu $, or
if $\mu ',\lambda $ do not commute, then by virtue of $\mu '$.

In the second case, 
note that $\lambda \in I(a_1,a_2)$ implies either
(a) $\nu $ is in the non-essential set of some $I(a_1',a_2')$
above $\mu $, (b) $\nu $ is separated from the lowest such 
$I(a_1',a_2')$ by a label with which it does not commute,
or (c) $\nu $ does not commute with some label
in the interior of the lowest such $I(a_1',a_2')$;
this follows from $\mu < \nu < \lambda < a_2$
together with the fact that $\mu $ and the label immediately above it
must form either a descent or an ILT.
In case 2(a), $\nu $ provides an alternative label for matching, and
$\lambda $ may be excluded from $S(a_1,a_2)$.  Notice that the partner 
cell in which $\nu $ or a preferable label $\nu '$
has been shifted may also exclude 
$\lambda $ by the following reasoning.  Since $\lambda $ commutes with 
$\nu $, the critical cell with $\nu\in I(a_1',a_2')$ and $\lambda $
immediately above $\mu $ may 
be matched by shifting $\nu $ downward to just
above $\lambda $, instead of by shifting $\lambda $.
  
For 2(b), include $\lambda $ in $S(a_1,a_2)$, since the cell
with $\nu\in I(\mu ,\lambda )$ may be matched by
shifting $\nu $ upward to its topologically decreasing
position above $I(\mu ,\lambda )$, or by shifting a label that is 
preferable to $\nu $.  Notice that we this also deals with case 4, 
by considering it from a different viewpoint.  
\begin{defn}\label{blocking-defn}
\rm
{
When such a label $\nu $ is matched by such
upward-shifting, we say that $\nu $ {\it blocks} $a_2$ from belonging to a
higher non-essential set. 
}
\end{defn}
Furthermore, observe that
shifting $\nu $ upward, or shifting a preferable label, still 
gives a critical cell which excludes $\lambda $ from $S(a_1,a_2)$; this is
because
either $\nu\in I(\mu ,\lambda )$, or $\nu $ appears above $\lambda $
and would form an ascend with $\mu $ if $\lambda $ were shifted upward.
For 2(c), if $I(a_1',a_2')$ has non-empty interior, then this gives an
alternative to $\lambda $, allowing $\lambda $ to 
be excluded from $S(a_1,a_2)$.  When this alternative label is shifted to
outside $I(a_1',a_2')$, it is still preferable to $\lambda $, by the 
conventions from case 2(a).  If there are no interior labels in
$I(a_1',a_2')$, then $\nu $ does not
commute with $a_2'$; $\lambda $ is included in $S(a_1,a_2)$, 
noting that the cell with $\nu $ shifted into $I(\mu ,\lambda )$ will be
matched by shifting $\nu $ upward to just below $a_1'$, similarly to 
case 2(b).

Now we turn to the third case.  First notice that $\mu $ must appear in a
lower ILT, either as its lowest label, or in its proper interior. 
If $\mu $ appears in the interior of some $I(a_1',a_2')$,
then exclude $\lambda $ from $S(a_1,a_2)$ and
apply our argument to $\mu $ or a preferable label,
proceeding downward until we find a way of matching.  By virtue of 
case (1), the matching partner which has shifted $\mu $ or a preferable 
label will also have excluded $\lambda $ from $S(a_1,a_2)$.
If $\mu $ is the lowest label of some $I(\mu ,\mu ')$ with non-trivial
interior, then again use a label from the interior or a preferable label
for matching.  On the other hand, if $I(\mu ,\mu ')$ has no interior,
then consider the descending labels $\nu_j,\nu_{j-1},\dots ,\nu_1 $
immediately above $\mu '$, up through the lowest label $\nu_1$ in the
next lowest ILT above $I(\mu ,\mu ')$.  Let $T$ be the subset 
of $\{ \lambda ,\nu_2,\dots ,\nu_j,\mu '\}$ which consists of those $\mu '$
along with those labels which
commute with $\mu '$, but not with $\mu $.  If all labels in 
$T$ belong to the non-essential sets of higher ILTs, then there is a 
Boolean algebra of critical cells in which each $T'\subset T$
specifies which labels 
to leave in decreasing order immediately above $\mu $, rather than shifted
upward into the interior of various ILTs; however, the 
empty set is missing from this Boolean algebra 
unless $\mu $ forms a descent with the 
label just above it when all labels in $T$ are shifts upward into interiors
of ILTs.  

Thus, allowing $\lambda \in S(a_1,a_2)$ and matching by shifting $\lambda $
to just above $\mu $ gives nearly 
a complete matching on this
Boolean algebra of critical cells,
but there is no
matching partner for the cell indexed by $T' = \{\lambda \} $.
However, the cell indexed by $\{ \lambda \} $ may 
instead be matched based on any
$\nu_i\not\in T$, since such a $\nu_i$
either belongs to the non-essential set of some higher ILT, or may 
be shifted downward into $I(\mu ,\lambda )$; in the latter case, we are
in a situation where $\nu_i $ blocks $\lambda $ from belonging to 
$S(a_1,a_2)$, in the sense described above.  In any event, all cells in
question are matched, and it is clear that 
the matching partners are also matched in the same fashion.

The fourth case was already handled within the argument for the
second case.
\EOP

\subsection{The Cohen-Macaulay Property}\label{cm-subsection}

In this section we
verify that the matching of the previous section consists of
covering relations in $F(\Delta)^M$, and that these comprise an acyclic 
matching on $F(\Delta)^M$ with only top-dimensional surviving critical cells.  
We begin with an important special case which captures most of the idea.

\begin{defn}
\rm
{
The {\it expanding interval} of a saturated chain $C$, denoted
$I(a_1,a_2)$, is the highest syzygy interval with non-empty non-essential
set in $C$..
}
\end{defn}

\begin{thm}\label{square-free}
If $I_{\Lambda }$ has a quadratic Gr\"obner basis and each $\phi^{-1}(m)$ is
square-free, then the poset interval $(0,m)$ has a 
discrete Morse function whose critical cells are all saturated, implying
the interval is homotopy equivalent to a wedge of spheres.
\end{thm}

\proof
In a syzygy interval labeled $a_1n_1\dots n_ja_2$,
recall that $a_1a_2$ is a Gr\"obner basis leading 
term.  Any unsaturated critical cell 
has at least one syzygy interval
with $j>0$.  We described in the previous section how to match all such 
cells so that partner cells differ by exactly one in their 
number of minimal skipped intervals.  
Lemma ~\ref{overlap} verifies that they in fact differ in dimension by
exactly one.  See Figure ~\ref{des_shift} 
for an example of a gradient
path from one such critical cell to its matching partner.

Lemma ~\ref{bool} shows for the expanding interval $I(a_1,a_2)$ that
$S(a_1,a_2)$ gives
rise to a Boolean algebra of critical cells within $F(\Delta)^M$, indexed 
by the subsets of $S(a_1,a_2)$.
$S(a_1,a_2)$ is chosen so that each $T\subseteq
S(a_1,a_2)$ gives rise to a unique such cell, denoted $Crit(T)$.
$Crit(T)$ is contributed by a saturated chain $M(T)$,
which has exactly the labels in $T$ inside 
$I(a_1,a_2)$, and each of the labels in $S(a_1,a_2) \setminus
T$ shifted to its topologically decreasing position outside $I(a_1,a_2)$.  
All labels in $M(T)$ 
other than $a_1,a_2$ and the members of $S(a_1,a_2)$ will appear in 
the same relative order for all choices of $T\subseteq S(a_1,a_2)$.

Any Boolean algebra has a complete acyclic matching simply by matching
by including/excluding any fixed set element.
We assign each critical cell to the Boolean algebra given by the 
non-essential set of its expanding interval, then take a union of 
complete acyclic matchings on these Boolean algebras. 
Lemma ~\ref{bool_group} checks that when one critical cell is assigned
to a particular Boolean algebra, then all critical cells in that Boolean
algebra are assigned to it, ensuring the matching is well-defined.  
Section ~\ref{match-subsection} 
already showed that we match all unsaturated critical cells.

The final step is to show
that this union of complete matchings on Boolean algebras
is an acyclic matching on $F(\Delta)^M$.  By Theorem ~\ref{multi-graph}, this would 
imply that we may simultaneously reverse all these gradient paths to
cancel all but some top-dimensional critical cells.  
To get acyclicity, we show two things: (1) Lemma ~\ref{filter} ensures that 
cycles cannot involve Boolean algebras from distinct fibres 
$\phi^{-1}(m)$, due to the filtration
$m_1\subseteq m_1\cup m_2 \subseteq \dots$ based on the monomial
term order $m_1,m_2,\dots $, and (2) Lemma ~\ref{boolcycle} verifies that
cycles cannot
involve multiple Boolean algebras in the same fibre.  Thus, we will 
produce a discrete Morse function whose critical cells are all 
top-dimensional, implying the order 
complex has the homotopy type of a wedge of spheres of top dimension.
\EOP

Next we deal with the possibility that not all monomials are
square-free.  The lemmas that follow do not use the square-free
assumption, so they apply to the general case.

\begin{thm}\label{morse-final}
If $I_{\Lambda }$ has a quadratic Gr\"obner basis, then the 
monoid poset $\Lambda $ has a discrete Morse function whose 
critical cells are all top-dimensional, implying $\Lambda $
is homotopically Cohen-Macaulay.
\end{thm}

\proof
The only issue left to address is repetition of labels.
To this end, we adjust the 
definition of non-essential set and make sure surviving critical cells 
still do not have any syzygy intervals with non-empty interior.  
When multiple copies of
a letter appear inside $I(a_1,a_2)$ or are upward-shiftable into it, 
only include one copy in
$S(a_1,a_2)$ that shifts downward to below $I(a_1,a_2)$;
we cannot shift more than one copy outside the interval
and still get a critical cell, since consecutive identical labels not
within a syzygy interval give a saturated chain rank not covered by
any minimal skipped interval.  Including one copy of the repeated letter
in the non-essential set is enough to ensure the Boolean algebra $B_n$
has $n\ge 1$, 
hence has a complete matching.
A letter cannot initiate or conclude a syzygy
interval and also appear in its interior, since then the syzygy interval
would not be a minimal skipped interval.  
We may have a syzygy interval which begins and 
ends with the same label $a_1$, but then there cannot be any interior 
labels at all.  
\EOP

\begin{lem}\label{overlap}
Critical cells that are matched
differ in dimension by exactly one.
\end{lem}

\proof
If none of the $I$-intervals are discarded in their conversion to 
$J$-intervals, then there is no issue (see Section 
~\ref{poset-disc-review} for definitions).  
When there is discardment,
this means there are three or more overlapping $I$-intervals
such that a middle one is unnecessary for covering all ranks 
by $I$-intervals, so Gr\"obner
basis leading term elements for these intermediate $I$-intervals will each
belong to the non-essential set of a higher $I$-interval, ensuring
matching by shifting such an individual label to outside the collection
of overlapping ILTs.  This matching operation preserves
the number of $J$-intervals from ILTs
and alters by exactly one the number of $J$-intervals coming
from descents.  Thus, dimension changes by exactly one.
\EOP

\begin{lem}\label{bool}
The critical cells indexed by subsets of
$S(a_1,a_2)$, for $I(a_1,a_2)$ the 
expanding interval of a saturated chain $M(T)$,
have the same incidences in $F(\Delta)^M$ as a
Boolean algebra of subsets of $S(a_1,a_2)$.  That is, there is a unique
gradient path from $Crit(T)$ to $Crit(T\cup \{ i\})$ for each $T\subseteq
S(a_1,a_2)$ and each $i\in S(a_1,a_2)
\setminus T$, and these are the only gradient paths among critical cells
in $B_{S(a_1,a_2)}$.
\end{lem}

\proof
Since all saturated chains in a fibre have equal content, and our
facet order is content-lex, downward steps in a gradient path must
sort labels.  The critical cell $Crit(T)$ is obtained by arranging
labels in $T$ in increasing order within $I(a_1,a_2)$, and labels 
in $S\setminus T$ in unique topologically 
decreasing positions below $I(a_1,a_2)$ (or above $I(a_1,a_2)$, in the
special circumstance that $a_2$ must be ``blocked'' from upward-shifting
into another non-essential set), from which they may shift into 
$I(a_1,a_2)$.   First we exhibit for each pair
$T'=T\cup \{ i\} $ with $i\in S\setminus T$, that 
there is a gradient path from $Crit(T)$ to $Crit(T')$.  
Choose $r$ and $t$ so that the $t$-th element in the chain $Crit(T)$ is 
just above the lowest label of
$I(a_1,a_2)$, and the $r$-th element of 
$Crit(T)$ is just above the label $i$ in $Crit(T)$, for $i$ satisfying 
$T'=T\cup \{ i\} $. 
In Example ~\ref{irt-examp}, let $i=3,t=3$ and $r=1$.

\begin{examp}\label{irt-examp}
Figure ~\ref{des_shift} depicts a gradient
path from a critical 2-cell of rank set $\{ 1,2,3\} $ to a 
critical 1-cell with rank set $\{ 1,2\}$ , based on a ring with 
$z_1z_4\in \init (I_{\Lambda })$.  We have $S=\{ 2,3\} $,
$T=\emptyset$ and $T'=\{ 3\} $.
\end{examp}

There is a gradient path 
from $Crit(T)$ to $Crit(T')$ of the form
$$d_r\circ u_r\circ d_{r+1}\circ
u_{r+1}\circ\cdots\circ d_{t-1}\circ u_{t-1} \circ d_{t},$$
because $i$ commutes with all labels separating it
from $I(a_1,a_2)$, $i$ is larger than all these separating
labels, and our discrete Morse function comes from a least-increasing
facet order.
Since the resulting permutation on labels is 321-avoiding,
Theorem ~\ref{one-path} ensures that this gradient path is unique, 
whether or not non-saturated chain segments are encountered in it.

To show that there are no other covering relations in $F(\Delta)^M$, i.e.
none between other pairs of
critical cells corresponding to subsets of $S(a_1,a_2)$, 
we use the fact that gradient paths can never introduce inversions.
Thus, a gradient path from $Crit(T)$ to 
$Crit(T')$ would imply $T\subseteq T'$, since any
$j\in T\setminus T'$, would imply an inversion $(j,a_1)$ in $Crit(T)$ 
that is not present in $Crit(T')$.
\EOP

Next we verify that critical cells are indeed partitioned into
Boolean algebras.

\begin{rk}
Critical cells with
label sequences of distinct content or with non-commuting labels in
opposite order are assigned to distinct Boolean algebras.
\end{rk}

\begin{lem}\label{bool_group}
Whenever one critical cell is assigned to a Boolean algebra, then all 
critical cells in that Boolean algebra are assigned to it.
\end{lem}

\proof
We must show 
that if $I(a_1,a_2)$ is the expanding interval for 
a saturated chain $M(T)$
for some $T\subseteq S(a_1,a_2)$, then the saturated chain
$M(T')$ for each $Crit(T')$ with $T'\subseteq S(a_1,a_2)$
also has $I(a_1,a_2)$ as its expanding interval.
Each label in a saturated chain
belongs to the non-essential set of at most one
syzygy interval, since it cannot pass through the lowest such 
syzygy interval above it to reach higher ones via a gradient path;
a label is only 
assigned to the non-essential set of a 
syzygy interval below it when it cannot shift into one above
it.  Shifting labels belonging to $S(a_1,a_2)$ from 
within $I(a_1,a_2)$ to their topologically
decreasing positions or vice versa
cannot cause a higher non-essential set to become
non-empty, by virtue of the choices made in Theorem ~\ref{well-def-bool}.
\EOP

\medskip
Finally, let us confirm that these complete 
matchings on Boolean algebras collectively
give an acyclic matching on $F(\Delta)^M$.  

\begin{lem}\label{boolcycle}
The matching on critical cells in $F(\Delta)^M$ is acyclic.
\end{lem}

\proof
Lemma ~\ref{filter} ensures there are no
directed cycles involving multiple fibres.  Suppose there were a 
directed cycle $C$ in a single fibre.  Any such $C$ must alternate upward
(matching) steps with downward steps.  Our matching consists of
a union of complete matchings on Boolean algebras.  Since the upward
steps in a fixed Boolean
algebra all insert the same fixed element $i$,
each downward step must take us to a different Boolean algebra, to 
avoid yielding the top of an upward-oriented edge in the same Boolean
algebra, from which the cycle could not have continued.

Suppose a matching step in $C$ shifts a 
label $\mu $ upward from within an ILT $I(a_1,a_2)$
to above it.  Then by virtue of our matching, $a_2$ must belong to the 
non-essential set of a higher ILT in the cell which has $\mu $ and 
all other labels within $I(a_1,a_2)$ shifted to below 
$I(a_1,a_2)$.  Furthermore, $\mu $ must not be upward-shiftable into
a higher ILT.  To pass to a distinct Boolean algebra, the 
downward step immediately after this upward-shifting of $\mu $
must either (1) shift $a_2$ upward into a higher ILT, (2) shift
a label $\lambda $ downward from a topologically decreasing
position 
into the interior
of an ILT, or (3) shift a label $\mu '$ upward 
into the interior of an ILT it then blocks 
(cf. Definition ~\ref{blocking-defn}).
(1) is impossible because $a_1,\mu $ form a non-inversion
after $a_2$ is shifted upward, implying an ascend between consecutive
commuting labels somewhere between $a_1$ and $\mu $, causing the cell not
to be critical.
In case (3), we can never un-do this
shifting of $\mu '$, since a matching
step will not shift it downward from $I(a_1,a_2)$, since 
$\mu $ blocks $a_2$ from shifting upward,
but $\mu $ is too large to shift below $a_2$ without $a_1$ also present.
Case (2) is allowed, but eventually we still would need to shift $\mu $
downward into $I(a_1,a_2)$, at which point we would have a downward step
keeping us in the same Boolean algebra, making it impossible for the 
cycle to continue.
Thus, we can rule out upward-shifting gradient path steps in a cycle.

Next suppose there is a step that shift labels downward either creating 
or eliminating an ILT.  Consider the lowest ILT $I(a_1,a_2)$ ever 
created/destroyed.  It must be destroyed by a downward step shifting
$a_1$ into the interior of a lower ILT $I(b_1,b_2)$.  Eventually we have
an upward (matching) step, shifting $a_1$ back upward from within 
$I(b_1,b_2)$ to below $a_2$.  But we have already eliminated the 
possibility of such upward-shifting matching steps within a cycle.

Finally, if all upward (matching) steps 
shift labels downward from within  
ILTs to between them, and all downward (non-matching) steps 
shift labels downward from between ILTs into lower ILTs, preserving 
the set of ILTs at each step, then labels not initiating or concluding
ILTs move progressively downward and may never return upward, making 
completing a cycle impossible.  To be precise we create inversions
between labels initiating/concluding ILTs and other labels, but we may
never eliminate these inversions.
\EOP

\section{Applications: minimal free cellular resolution and a finite state
automaton which computes Poincare'-Betti series}\label{section-quad-applic}

This section describes the surviving critical cells in the quadratic
Gr\"obner basis case in two ways:
\begin{enumerate}
\item
as the words generated by a finite state automaton, 
implying the generating function for Morse numbers
is rational
\item
as representatives of 
the $J'$-non-stuttering, $J'$-commuting equivalence classes of words,
as developed in [HRW] to count Betti numbers.
\end{enumerate}

\begin{thm}
The discrete Morse function of the previous section
gives a minimal free cellular resolution of $k$ as a 
$k[\Lambda ]$-module.
\end{thm}

\proof
Results of [BW] imply that the complex of critical cells from our
Morse function supports a free cellular resolution of $k$ as a 
$k[\Lambda ]$-module, because our acyclic
matching preserves multi-grading.  Furthermore, there are no 
incidences among critical cells of equal multi-degree
in this Morse function, because all 
critical cells of multidegree $\lambda $
come from saturated chains with highest element $\lambda $, 
making gradient paths from
one critical cell to another of the same multi-degree
impossible, despite the fact that critical 
cells need not be concentrated in a single dimension.  
If we consider 
the complex obtained by tensoring the complex
of critical cells with $k$, this implies that all its boundary
maps are 0 maps.  This implies that the resolution supported by
the complex of critical cells is a minimal free resolution.

Alternatively, one may see that the resolution is minimal by checking
that Morse numbers equal Betti numbers.
Theorem ~\ref{hrw-biject} does this by constructing a 
bijection between the critical cells in our Morse function and the 
$J'$-non-stuttering, $J'$-commuting equivalence classes of words of [HRW];
these equivalence classes of words were shown in [HRW] to 
index a basis for $\rm{Tor}(k,k)$, because those of fixed multidegree
$\lambda $ index a homology basis for $(0,\lambda )$.

\EOP

\begin{rk}\label{morse-betti-rk}
The fact that the discrete Morse function gives a minimal free cellular
resolution implies Morse numbers equal Betti numbers.
Thus, the upcoming generating
function for Morse numbers also computes the Poincare'-Betti  series.  
\end{rk}

The next theorem 
will construct a finite state automaton that generates exactly the 
label sequences for the surviving critical cells.
The list of states in this finite state automaton is 
far from minimal in general among all finite state automata
generating this language.  Specifically, we keep 
track of more data in each state than is strictly necessary, in
order to greatly simplify the description of our automaton.

\begin{thm}\label{finite-state-auto}
The label sequences for saturated chains which contribute
surviving critical cells are exactly the words of a regular language.
Thus, the generating function for Morse numbers, which in this 
case equals the Poincare'-Betti series, is rational.
\end{thm}

\proof
The alphabet for the language is the set of labels on covering relations,
i.e. of generators for the monoid.  
For convenience, we view label sequences on saturated chains
as words by reading them from top to bottom.  Since all 
surviving critical cells are saturated chains, the dimension of
each such critical cell is two less than the length of the word labelling 
it.  Thus, the Morse number $m_i$ counts words of length $i+2$ in the 
language of label sequences.  
We will describe a set of states
and of legal transitions between states that comprise
a finite state automaton that generates exactly the language of label
sequences for surviving critical cells.
The existence of such an automaton will imply that the 
language is regular, and hence the generating 
function for Morse numbers is rational (cf. 
Section ~\ref{regular-review}).  In fact, the rational generating 
function may be determined from the finite state automaton (see [BR]).
The remainder of the proof describes how to construct such an automaton.

The automaton has a 
unique initial state, and each time a label 
is read, a transition is 
made from one state to another state if the label sequence read so far
could be the initial segment for a label sequence of a surviving critical
cell.  To decide which labels give valid transitions, each state 
must keep track of enough data about previously read labels to decide whether
concatenating a newly read label will
\begin{itemize}
\item
give a label sequence for a surviving
critical cell, in which case a transition is made to a final state, or
\item
give a label sequence for a critical cell which is cancelled, but one
where reading additional labels could again yield a surviving critical
cell, in which case a transition is made to a non-final state, or
\item
give a label sequence not meeting either of the above forms, in which
case there is no valid transition, so the word is not generated by the 
automaton.
\end{itemize}

Specifically, for there to be a transition labeled $\lambda $ out 
of a state $S$, $\lambda $ 
must form either a descent or an ILT with the most recently read label,
and there are further constraints related to non-essential sets.
The requirement about descents and ILTs is necessary because
every pair of consecutive labels for a surviving critical cell must take
this form.  

Each state will contain the following data:
the list of previously encountered ILTs and individual 
labels, together with the order of the most recent occurrences of these
ILTs and labels.  Thus, each state has associated to it 
a subset of the finite set of monoid generators and
leading terms in our Gr\"obner
basis, together with a permutation on the elements of this subset.
Earlier occurrences of the same ILTs or 
individual labels are unnecessary for deciding whether all non-essential 
sets are empty, or else would have already caused the word to be 
unproducable by the automaton at an earlier stage.  Thus, we have a 
finite list of states.

If a partial label sequence $w $ concludes with a label $\mu $
and leads to a final state $S$, then 
the next label $\lambda $ to be read gives a legal 
transition from $S$ to another final
state if and only the following conditions are all met:
\begin{enumerate}
\item
$\lambda ,\mu $ comprise a descent or ILT 
\item
$\lambda $ is not in the non-essential set of any earlier ILT.  That is,
every previously encountered 
ILT $I(a_1,a_2)$ with $a_1 < \lambda <a_2$ either has
(a) $\lambda \nu \in \init (I_{\Lambda })$ for some label $\nu $ 
read more recently than $I(a_1,a_2)$, (b) $\lambda a_1\in \init 
(I_{\Lambda })$, (c) $\lambda a_2\in \init (I_{\Lambda })$, or
(d) $\lambda $ smaller than some label $\nu $ read after
$I(a_1,a_2) $
\item
If $\lambda ,\mu $ comprise an ILT, then there is no previously
encountered label $\mu '$ in its non-essential set.
That is, there is no previously encountered 
$\mu '$ 
with all the following properties: 
(a) $\mu '$ satisfies $\lambda < \mu ' < \mu $, 
(b) $\mu '$ is smaller than all labels read after it and before $\lambda $, 
(c) $\mu '$ commutes with all labels read after it,
and (d) deleting $\mu '$ would cause $\mu $
to be in the non-essential set of a previously encountered ILT.
\end{enumerate}

When the first and third conditions hold but the second one
fails, there is still a transition to a 
non-final state $U$. 
However, the only legal transitions from such a non-final state $U$ 
are given by labels
$\lambda '$ such that $\lambda ',\lambda $ form an ILT which causes
$\lambda $ no longer to belong to a non-essential set, i.e. 
when we are in one of the following circumstances:
\begin{itemize}
\item 
$\lambda ' < \mu $ and $\lambda '\mu\not\in \init (I_{
\Lambda })$, because then shifting $\lambda $ upward would yield a 
non-critical cell
\item
$\lambda '$ would also belong to the non-essential set of some 
ILT once $\lambda $ is shifted upward into an ILT (i.e. 
Theorem ~\ref{well-def-bool}, case 1), implying the critical cell with
ILT $\lambda ',\lambda $ is not matched by shifting $\lambda $ 
upward.
\end{itemize}

The necessity of these constraints on allowable
words is immediate from the description of critical cells and the 
matching to cancel them in earlier sections.  
These constraints on legal transitions 
are also sufficient to produce a surviving
critical cell because any such label sequence will label a critical
cell whose ILTs all have empty non-essential set, i.e. a critical cell 
that is not cancelled.
\EOP

\medskip
Following [HRW], let $J=\init (I_{\Lambda })$, and 
let $J'$ be the complement of $J$.  
Labels $a,b$ commute if and only if $ab\not \in J$, in which case we say
they are {\it $J'$-commuting}.  Define a 
{\it $J'$-commuting equivalence class} of label sequences to be a set 
of label sequences which agree up to $J'$-commutation.
A label sequence is
{\it $J'$-stuttering} if it has consecutive labels
$a_1,a_1$ where $a_1^2\not\in J$.  
A $J'$-commuting equivalence class $C$ is {\it $J'$-non-stuttering} if none of
the label sequences in $C$ are $J'$-stuttering.

\begin{thm}\label{hrw-biject}
There is a bijection between the 
$J'$-non-stuttering, $J'$-commuting equivalence classes of a given 
content and the label sequences of the same content 
for critical cells that survive 
cancellation.  Moreover, exactly one member of
each $J'$-non-stuttering, $J'$-commuting equivalence class is a 
label sequence for a critical cell surviving cancellation.
\end{thm}

\proof
We will show that each $J'$-non-stuttering $J'$-commuting equivalence class
contains exactly one label sequence for a critical cell surviving 
cancellation.  First we show the existence of such a label sequence 
within each such $J'$-non-stuttering $J'$-commuting
equivalence class by
providing an algorithm which applies a series of $J'$-commutation 
relations to transform any member of such a class into 
the label sequence for a 
surviving critical cell.  Then we show that each such class has at most
one label sequence from a critical cell surviving cancellation.
Finally, we show that $J'$-stuttering anywhere within a 
$J'$-commuting equivalence class of a label sequence implies that
the label sequence either does not come from a critical cell or is
cancelled.

The algorithm sequentially processes the labels, proceeding from 
smallest to largest label value, and in the case
of repetition, proceeds from highest to lowest initial location for
each value.  The algorithm terminates
because it processes a finite number of labels and will 
use a finite number of steps to process each label.
If the label $\mu $ immediately below a label $\nu $ to be processed is 
smaller than $\nu $, then the algorithm would have processed
$\mu $ before $\nu $, and we will soon see that 
the pair must form an ILT.  In this case, we say 
that $\mu $ is {\it attached} to $\nu $ at the time $\nu $ is processed.  

A label $\nu $ is processed as follows.  If $\nu $ is not attached to a 
label immediately below it, 
then $\nu $ is shifted upward until $\nu $
either encounters a label with which it does not commute, reaches 
the top of the label sequence, or encounters a label smaller than it 
such that all current ILTs $I(a_1,a_2)$ above $\nu $ with $a_1 < \nu < a_2$
and $a_1\nu ,\nu a_2\not\in \init (I_{\Lambda })$
are currently separated from $\nu $ by labels with which $\nu $ does not
commute.
Notice that $\nu $ will not encounter another copy of $\nu $ before
reaching such a position, because the $J'$-commuting equivalence class is
$J'$-non-stuttering.
If $\nu $ is shifted to just below a label $\lambda_2$ with which 
$\nu $ does not commute, such that $\nu < \lambda_2$, then $\nu $ is now
attached to $\lambda_2$.
It in addition the label $\lambda_1$ previously below
$\lambda_2$ had formed an ILT with $\lambda_2$, then $\nu < \lambda_1 <
\lambda_2$, and we may detach $\lambda_1$ from $\lambda_2$ at the same
time that we attach $\nu $ to $\lambda_2$ to form a new 
ILT $I(\nu ,\lambda_2)$.
If a label $\nu $ to be processed is attached to 
a label $\mu $ just below it,
then the pair is shifted upward as a unit past labels larger than $\nu $
that commute with both $\mu $ and $\nu$, with
the following special rules: 
\begin{itemize}
\item
if $\mu ,\nu $ encounter a label
$\lambda > \nu $ which commutes with $\nu $ but not with $\mu $,
then detach $\nu $ from $\mu $, shift $\nu $ past $\lambda $, 
and attach $\lambda $ to $\mu $, and continue processing $\nu $ as an
unattached label
\item
if the label immediately above $\nu $ forms either a descent or ILT
with $\mu $ and $\nu $ can be shifted upward into the interior of and 
ILT $I(a_1,a_2)$ with $a_1<\nu <a_2$, $a_1\nu,\nu a_2\not\in \init 
(I_{\Lambda })$, and $\nu $ commuting with all labels between $\nu $ and
$I(a_1,a_2)$, then $\nu $ is detached from $\mu $ and shifted upward to
above $I(a_1,a_2)$ and continues its processing.
\end{itemize}

The fact that the only ascends in the output are between non-commuting pairs
ensures it labels a critical 
cell.  To see it also is one that survives cancellation, one may check 
that all non-essential sets are empty.  At the time $\mu $ was 
processed, $\mu $ could not be further shifted upward into an ILT for
which it would belong to the non-essential set, or else $\mu $ would have
been shifted farther upward in its processing.
The fact that all smaller values had already been processed by the time
$\mu $ was processed 
ensures that this property is preserved throughout the algorithm.
We may also eliminate the possibility of non-essential set members
that shift downward into ILTs, because these only
arise when the top of some ILT is capable of shifting upward without its
partner, but our algorithm would have actually performed this
shifting, and again the fact that we process smaller labels before 
larger ones means this property is preserved throughout the rest of 
the algorithm.  Thus, all non-essential sets are empty, so the algorithm
indeed outputs the label sequence of a surviving critical cell.  

To show that there is at most one label sequence surviving 
critical cell cancellation in each $J'$-commuting equivalence class,
first note that pairs of consecutive labels that commute must appear
in descending order to avoid either
having an ascend not appearing within an ILT
(implying the saturated chain does not contribute a critical cell)
or having an ILT with non-empty interior (implying the critical cell is
cancelled).  Now suppose there are two 
label sequences in the same $J'$-commuting equivalence class, both
from critical cells surviving cancellation.  Then there
must be some pair of $J'$-commuting labels $\mu_1,\mu_2$ with $\mu_1
< \mu_2$, such that the pair  
are inverted in one label sequence and not the other.  
One may use the intermediate
value theorem to show that the non-inverted pair $\mu_1 ,
\mu_2$ must be separated by at least one ILT $I(a_1,a_2)$ with 
$a_1 < \mu_1 < a_2$ and by at least one ILT $I(b_1,b_2)$ with 
$b_1 < \mu_2 < b_2$, such that $I(a_1,a_2)$ either equals $I(b_1,b_2)$
or occurs before $I(b_1,b_2)$.
This implies either that $\mu_1\in S(a_1,a_2)$, ensuring cancellation, or 
that some label not commuting with $\mu_1$ separates 
$\mu_1$ from $I(a_1,a_2)$, 
implying $\mu_2$ must commute with all labels between it and $a_1$, 
in order for it to be possible to swap $\mu_1,\mu_2$.
But then the critical cell with $\mu_2, \mu_1$ forming an inversion
must have $\mu_2 \in S(b_1,b_2)$, and so cannot also survive cancellation.

Now we turn to the issue of non-stuttering.
Any two consecutive identical labels must appear within 
a syzygy interval to avoid comprising an ascend which would make the cell
non-critical.  But we showed that any critical cell with a syzygy
interval with non-empty interior is cancelled. 
If a label sequence $\omega $ is 
in the $J'$-commuting equivalence class of a 
label sequence which has consecutive identical labels, 
then some label $\lambda $ appears more 
than once in $\omega $, separated by one or more syzygy 
intervals.  But then one
of these syzygy intervals will have the lower copy of $\lambda $ 
in its non-essential set, 
implying the critical cell is cancelled, unless $\lambda $ does not 
commute with some separating label.  But this label also does not 
commute with the other copy of $\lambda $, making it impossible for 
the two labels to be shifted to consecutive positions, as needed for
stuttering, a contradiction.  Thus, if a label sequence is 
$J'$-commuting
equivalent to one with stuttering, then the label sequence does
not survive cancellation.
\EOP

\section{Gr\"obner bases of higher degree}\label{high_deg}

In this section, we extend results of Section ~\ref{cm-section} from 
the quadratic Gr\"obner basis case to the degree $d$ to prove:

\begin{thm}\label{deg-d}
If $I_{\Lambda }$ has a Gr\"obner basis of 
degree $d$, then $\tilde{H}_i (\Delta (\hat{0},\lambda )) = 0 $ for 
$i < -1 + \frac{\deg (\lambda )-1}{d-1} $, with $\deg (\lambda )$ defined
as below.  Hence, $Tor_i^{k[\Lambda ]} 
(k,k)_{\lambda } = 0$ for 
$i< 1 + \frac{\deg (\lambda )-1}{d-1} $.  Moreover, this vanishing is
achieved by a free cellular resolution resulting from a discrete Morse
function on $\Lambda $.
\end{thm}

In the standard-graded case, $\deg (\lambda )$ is given by the 
grading.  In general, let $\deg (\lambda )$ be
one more than the length of the shortest 
saturated chain on the poset interval $(\hat{0},\lambda )$, i.e. the
degree of the image of $\underline{x}^{\lambda }$ in 
the associated graded ring.  

\begin{prop}\label{deg-d-start}
Ordering saturated chains by using any monomial
term order to order fibres then lexicographically ordering 
saturated chains
within each fibre yields a content-lex facet order.
\end{prop}

\proof
Syzygy intervals now must be defined to have weakly increasing labels 
$a_1,\dots ,a_r$ such that there is a Gr\"obner basis leading 
term which divides $a_1\cdots a_r$ and has smallest divisor $a_1$ and 
largest divisor $a_r$.   To be a minimal skipped interval, we must also
have that neither
$a_2\cdots a_r$ nor $a_1\cdots a_{r-1}$ is divisible by a Gr\"obner basis
leading term.  Then the proof of Theorem 
~\ref{morse-init} applies.
\EOP

The remainder of this section is concerned with
cancelling pairs of critical cells to obtain a Morse function
with no critical cells below dimension $-1 + \frac{\deg (\lambda )-1}{d-1}$.

\begin{defn}
\rm
{
An {\it increasing leading term}, or {\it ILT}, is a Gr\"obner basis leading
term, with labels arranged in weakly increasing order.  We will often use 
the term ILT to refer to an ILT that labels a specific syzygy interval.
}
\end{defn}

We will use the variables $d',d''$ to 
represent the degree of an arbitrary Gr\"obner basis leading term, so we
always will have $d',d''\le d$.
Denote a syzygy interval with ILT $a_1\cdots a_{d'}$ by
$I(a_1,\dots ,a_{d'})$.
In contrast to the quadratic Gr\"obner basis case, now 
there may be several Gr\"obner
basis leading terms specifying the same syzygy interval.
This fact, that a single syzygy 
interval may have several ILTs beginning and ending with the same pair of
labels, causes one substantial new issue to arise: the critical cells
resulting from one syzygy interval may comprise
several overlapping Boolean algebras,
since
different ILTs will give rise to different non-essential sets, as in
Example ~\ref{overlap-bool-examp}.
\begin{examp}\label{overlap-bool-examp}
Consider the syzygy interval labeled $a_1,\dots ,a_4$ with 
Gr\"obner basis leading terms $a_1a_2a_4$ and $a_1a_3a_4$.
One Boolean algebra of critical cells, based on ILT $I(a_1,a_2,a_4)$
consists of labels sequences $a_1a_2a_3a_4$ and $a_3a_1a_2a_4$, while
the Boolean algebra for $I(a_1,a_3,a_4)$ consists of label sequences
$a_1a_2a_3a_4$ and $a_2a_1a_3a_4$.  The critical cell labeled 
$a_1a_2a_3a_4$ is shared by the two Boolean algebra.
\end{examp}
Lemma ~\ref{mult-ilt} deals with such overlap by providing an acyclic 
matching for any such collection of overlapping Boolean algebras.

Define the non-essential set for $I(a_1,\dots ,a_{d'})$, denoted
$S(a_1,\dots ,a_{d'})$, similarly to the $d=2$ case, but now it may have
two types of members:
\begin{enumerate}
\item
individual labels, in exact analogy 
to the non-essential set members for the $d=2$ case, i.e. labels
which either (a) appear in topologically decreasing positions below 
$I(a_1,\dots ,a_{d'})$ (or above $I(a_1,\dots ,a_{d'})$ in the 
exceptional circumstances discussed in Theorem ~\ref{well-def-bool})
from which they may shift into $I(a_1,\dots ,
a_{d'})$ via a gradient path without causing $I(a_1,\dots ,a_{d'})$ 
to cease to be a minimal skipped interval, or (b)
labels that have thus shifted into the interior of $I(a_1,\dots ,a_{d'})$.
\item
collections $\{ b_2,\dots ,b_{d'} \} $ of labels that appear 
either immediately above a label $b_1$ with which they form an ILT, or
which collectively appear
in the interior of one or more ILTs strictly above such a label
$b_1$, with
$I(a_1,\dots ,a_{d'})$ serving as the highest of these ILTs.
Furthermore, we require there to be a gradient path from the former to
the latter which shifts the collection of labels upward into the interior
of the various ILTs in order for the collection of labels to belong to
$S(a_1,\dots ,a_{d'})$.
\end{enumerate}

Follow the conventions of Theorem ~\ref{well-def-bool}
to decide which such individual labels and collections of labels 
should belong to $S(a_1,\dots ,a_{d'})$, not allowing 
collections of labels that shift downward into $I(a_1,\dots ,a_{d'})$ 
from above.  If there is a need for ``blocking'' as in Theorem
~\ref{well-def-bool}, an individual 
label will always serve this function rather than a collection of labels.  
If a label $\lambda \in I(a_1,\dots ,a_{d'})$ meets the above requirements
to be included individually in $S(a_1,\dots ,a_{d'})$, then it is not also 
included as part of a collection of labels.
With these conventions, the proof that this gives a 
well-defined matching is identical to the proof of
Theorem ~\ref{well-def-bool}.

\proofo
Proposition ~\ref{deg-d-start} 
provides the Morse function that serves as our starting point.
We will cancel all critical cells of dimension less than
$-1 + \frac{\deg (\lambda )-1}{d-1} $, using a fairly similar, but 
somewhat more subtle, approach
to the $d=2$ argument.  Notice that each critical cell of dimension
less than $-1 + \frac{\deg (\lambda )-1}{d-1} $ is contributed by
a saturated chain with average minimal skipped interval height
greater than $d-1$.   But any minimal skipped interval of
height greater than $d-1$ is a syzygy interval with more than $d$ 
labels, so it is consists of 
an ILT with at least one additional label
interspersed.  Such extra labels either allow the critical cell to be 
cancelled similarly to the $d=2$ case, or in the case of a
syzygy interval with multiple ILTs, Lemma ~\ref{mult-ilt} gives a 
matching in which all unmatched cells
have average interval height at most $d-1$ for the $I$-intervals
related to the syzygy interval.  Thus, cells left unmatched
must then have another syzygy interval of height greater than $d-1$
at lower ranks.  This allows us to repeat the argument until 
eventually reaching a syzygy interval which causes the cell to 
be cancelled.

The fact that there is indeed a unique gradient path from a critical
cell $\tau $ to a critical cell $\sigma $ for each pair $\tau ,\sigma $
to be cancelled follows from Theorem ~\ref{one-path}.  When a single 
label is shifted into a syzygy interval, the gradient path is 
identical to the one given in 
the $d=2$ case.  When a collection of labels is shifted upward into a 
syzygy interval, the gradient path is the one  
described in the proof of Theorem ~\ref{one-path}.  
Theorem ~\ref{one-path} also proves the uniqueness of these 
gradient paths, whether shifting a single label or a collection of 
labels.  As before, we give a complete acyclic matching on each Boolean
algebra of critical cells, as long as it is not part of a collection of 
overlapping Boolean algebras.  Similarly to the $d=2$ case, we match all 
cells in this Boolean algebra by including/excluding a single 
non-essential set member $\nu $ from the interior of the syzygy interval.
It is convenient to choose $\nu $ to be the individual label with highest
topologically decreasing position outside the syzygy interval, if there
is such a label, and otherwise to choose the collection of labels with
highest topologically decreasing position outside the syzygy interval. 
Lemma ~\ref{mult-ilt} provides the matching for collections of overlapping
Boolean algebras.

With these choices, acyclicity is similar to the $d=2$ case, since
Lemma ~\ref{mult-ilt} will
verify acyclicity of the matching on a single collection
of overlapping Boolean algebras resulting from several ILTs on a single
syzygy interval.  Applying results of [BW],
the desired resolution is immediate from this acyclic matching.
\EOP

Next we prove Lemma ~\ref{mult-ilt} in the
special case of degree $d=3$.  In this case we deduce a stronger 
result than for general $d$, but the proof is also much simpler than
in general, but gives the flavor of the upcoming proof for degree $d$.

\begin{defn}
\rm
{
A pair of ILTs $a_1\dots a_r$ and $b_1\dots b_s$ 
with $a_1\le b_1,a_r\le b_r$ are {\it concatenating} if
either (1) $a_1<b_1$ and $a_r=b_i$ for some $1\le i<s$, or 
(2) $a_r<b_s$ and $b_1=a_i$ for some $1<i\le r$.
}
\end{defn}

\begin{lem}\label{mult-ilt-deg3}
Let $I_{\Lambda }$ be a toric ideal with Gr\"obner basis 
with leading terms all of degree at most 3.
Let $I=I(a_1,\dots ,a_{d'})$ be a syzygy interval, given by one or
more ILTs, each of which gives rise to a Boolean algebra of critical
cells.
Then this collection of overlapping Boolean algebras 
has an acyclic 
matching which matches all critical cells with average interval height
at most 2 for $I$ together with any 
descents coming from labels shifted out of $I$.
\end{lem}

\proof
Order the ILTs $M_1,\dots ,M_k$.  If some $M_i$ has degree 2, then the 
Boolean algebra for each $M_j$ is contained in the Boolean algebra for
$M_i$, so we use the complete matching on a single Boolean algebra.
Otherwise, we have labels $\lambda_1,\dots ,\lambda_k$ such that for
$1\le i\le k$, $M_i = a_1\lambda_ia_{d'}$, for fixed initial and final 
labels $a_1,a_{d'}$.  Let $\Delta_i$ be the set of critical cells in the
Boolean algebra $B^{(i)}$ for $M_i$ which are not shared with any 
earlier Boolean algebra $B^{(i')}$ for $i'<i$.
Notice that $\Delta_i$ 
consists of exactly
those critical cells in $B^{(i)}$ which have $\lambda_1,\dots ,
\lambda_{i-1}$ all shifted to topologically decreasing positions
outside $I$.  Thus, $\Delta_i$ has the structure of a Boolean algebra,
resulting from all other labels in the non-essential set for
$I(a_1,\lambda_i,a_{d'})$, so this has a complete acyclic matching unless
this set is empty.  But when the set is empty, then $I$ consists of 
only the three labels $a_1,\lambda ,a_{d'}$, as well as labels essential
to concatenating ILTs, so matching is not 
necessary.  In the case of concatenating ILTs, the average interval 
height is still at most 2.
\EOP

\medskip
The situation gets much more complex when labels other than 
$a_1$ and $a_{d'}$ 
may divide more than one of the Gr\"obner basis leading terms specifying
ILTs on the syzygy interval.
 
\begin{lem}\label{mult-ilt}
Suppose that a 
single expanding interval $I$ has multiple ILTs.  Then the resulting 
collection of overlapping Boolean algebras has an acyclic matching such
that all unmatched cells have minimal skipped interval average height
at most $d-1$.
\end{lem}

\proof
Choose a total order $M_1,\dots ,M_k$ on the ILTs for $I(a_1,
\dots ,a_{d'})$.  Thus,
each $M_i$ is a Gr\"obner basis leading term with smallest divisor
$a_1$, largest divisor $a_{d'}$ and with divisors of intermediate value,
all of which appear as labels that can shift in/out of the syzygy 
interval.  Thus, any shifting of labels in/out of the syzygy interval
still gives a syzygy interval as long as at least one of these ILTs 
appears entirely within the syzygy interval.
Each ILT $I(M_i)$ has its own non-essential set, denoted $S(M_i)$, 
giving rise to its own Boolean
algebra of critical cells.  
Denote by $\Delta_r$ the collection of critical cells in the 
Boolean algebra given by $M_r$ which are not shared with any of the 
earlier Boolean algebras given by $M_1,\dots ,M_{r-1}$.  

We will provide an acyclic matching on each such $\Delta_r$.
Note that $\Delta_r$ consists of those subsets of $S(M_r)$ which
shift enough labels to outside the syzygy interval $I(M_r)$ so that
the label sequence on $I(M_r)$ is not divisible by any
of the monomials $M_1,\dots ,M_{r-1}$.  If there is any label in 
$S(M_r)$ that does
not divide any of the monomials $M_1,\dots ,M_{r-1}$, then we obtain a 
complete acyclic matching on $\Delta_r$ 
by including/excluding one such label in $I(M_r)$.
Next we consider the case where each member 
of $S(M_r)$ does divide some earlier $M_i$.

Fix an ordering $\lambda_1,\dots \lambda_l$ on the elements of 
$S(M_r)$.  It will be convenient in the next section
if we order them from highest to lowest
topologically decreasing position outside $I(M_r)$.
Now apply the following matching
procedure to each critical cell in $\Delta_r$:
\begin{enumerate}
\item 
match the cell by including/excluding $\lambda_1$ from $I(M_r)$ 
unless shifting $\lambda_1$ to inside $I(M_r)$
yields a cell in an earlier $\Delta_i$,
\item 
if the cell is not yet matched, then 
match by including/excluding $\lambda_2$ from $I(M_r)$, 
unless this yields a matching partner which was
already matched at the first step or which belongs to an 
earlier $\Delta_i$
\item 
continue inductively, matching the cell
by including/excluding $\lambda_i$ from $I(M_r)$ if the cell
was not already matched based on any of the labels
$\lambda_1,\dots ,\lambda_{i-1}$ and the partner cell based on
shifting $\lambda_i$ also does not belong to an earlier $\Delta_i$ 
and is not already matched based on any earlier label $\lambda_{i'}$
with $i'<i$.
\end{enumerate}

Notice that a cell cannot be matched based on the label
$\lambda_i$ if either (a) $\lambda_i\not\in I(M_r)$
and shifting $\lambda_i$ into $I(M_r)$ gives a cell in an earlier Boolean 
algebra, or (b) shifting $\lambda_i$ in or out 
of $I(M_r)$ gives a cell previously matched.  
Thus, any unmatched critical cell that has exactly the 
labels $M = \{ \mu_1,\dots ,\mu_k \}$ shifted to outside $I(M_r)$ 
will have the 
property that each such $\mu_i$ is necessary outside $I(M_r)$ either to 
avoid overlap with an earlier Boolean algebra or in order for some
$\nu \in I(M_r)$ coming earlier than $\mu_i$ also not to allow matching.
That is, in the latter case there must be some $\mu_{i'} < \nu < \mu_i$ with 
$\mu_{i'}\not\in I(M_r)$,
such that $\mu_{i'}$ ``covers'' multiple ILTs (see Definition 
~\ref{cover-defn}), some of which could also be 
covered by $\nu $, and the rest of which are also covered by $\mu_i$.
\begin{defn}\label{cover-defn}
\rm
{
A label $\mu $ {\it covers} an earlier ILT $M_j$ if $\mu $ divides the 
Gr\"obner basis leading term specifying $M_j$.
}
\end{defn}
Assign to each $\mu_i\in M$ either an ILT $M_{i'}$ which 
it exclusively covers, 
or an ILT $M_{i''}$ that it would exclusively cover if the 
earliest forbidden $\nu \in I(M_r)$ 
were shifted to outside $I(M_r)$, or which it would exclusively cover 
after some number of iterations of this reasoning, i.e. an ILT which 
makes it impossible to shift $\mu_i$ from outside $I(M_r)$ to inside
$I(M_r)$ as a matching step.  Call this ILT which is 
assigned to $\mu_i$ the {\it indexing ILT} of $\mu_i$.

If we can show that every label
in $S(M_r)$ belongs to one of the $k$ indexing ILTs, this will 
imply $|T|\le k\cdot (d-2)$, as desired.  
Suppose some $\lambda \in S(M_r)$ is not in any of the 
indexing ILTs, and choose the 
label $\lambda $ of this form which comes earliest in our ordering
on labels in $S(M_r)$.
First note that $\lambda\in I(M_r)$, since otherwise $\lambda $ would
belong to its own indexing ILT.  
We will show next that the cell with $\lambda $ shifted to
outside $I(M_r)$ is not matched based on a label of higher precedence 
than $\lambda $.  Since shifting $\lambda $ to outside $I(M_r)$ also 
cannot give a cell belonging to an earlier $\Delta_i$, we will be able
to conclude that the critical cell will be matched based on $\lambda $.
Thus, any unmatched cell will satisfy $|T| \le k\cdot (d-2)$.

Now we prove the claim that the cell may be matched by shifting
$\lambda $ to outside $I(M_r)$.  
When $\lambda $ is shifted to outside $I(M_r)$,
each $\mu_i$ of higher precedence which appears
outside $I(M_r)$ cannot be shifted to inside $I(M_r)$ as a matching
step, by virtue of  
its indexing ILT, since $\lambda $ cannot cover this indexing ILT. 
Likewise any $\lambda_i$ of higher precedence
which appears within $I(M_r)$ in the critical 
cell cannot be shifted to outside $I(M_r)$ 
without rendering some $\mu_{i'}$ 
unnecessary for covering its indexing ILT, since otherwise we would 
have matched based on the smallest $\lambda_i$ which did not have
this property; in particular, this means that
$\lambda_i$ must belong to the indexing ILT for 
$\mu_{i'}$ in the critical cell.  Shifting $\lambda $ to outside $I(M_r)$
does not change this relationship, so the cell with $\lambda $ shifted
to outside $I(M_r)$ also cannot match by shifting $\lambda_i $.
Thus, $\lambda $ is the first label allowing matching for both cells, so
both are indeed matched by shifting $\lambda $.

In the case where all elements of 
$S(M_r)$ are individual labels that shift to topologically
decreasing positions outside $I(M_r)$, 
this yields the following upper
bound on average interval height for this portion of the interval
system, using the fact
that total height is one less than the total number of
labels involved, and that $k$ is non-negative:
$$ \frac{\rm{total \hspace{.1in} height}}{\rm{no. \hspace{.1in} of 
\hspace{.1in} intervals}} =
\frac{d-1 + |T|}{k+1}\le \frac{d-1 + k\cdot (d-2)}{k+1 } 
= (d-2) + \frac{1}{k+1}\le d-1.$$
Let us now handle the more general case, where
some non-essential set members are collections of
labels.  All labels belonging to such collections 
will contribute
individually to the bound $|T|\le k\cdot (d-2)$ when the labels
appear within $I(M_r)$, because each label
contributes individually to monomial degree.
When such a collection of labels appears outside the ILT, it 
would increase the number of intervals $k$
by one, but would increase the total height by as much as $d-1$,
seemingly invalidating the above computation of average interval
height.  However, the highest label in the
newly created ILT must also form a descent with the 
label immediately above it, and we may use this descent rather
than the new ILT in order to compute the above bound, since the descent 
will not also be counted in
a similar computation for any other syzygy interval.  We may
safely ignore the newly created ILT in the bound 
computation, since it also has height at most $d-1$.

Acyclicity will follow from the Cluster Lemma of 
[Jo] (see Lemma ~\ref{filter}), using the 
filtration of subcomplexes
$B_1\subseteq B_1\cup B_2\subseteq \cdots \subseteq
B_1\cup \cdots \cup B_k$ where $B_1\cup \cdots \cup B_i$ is the union of
Boolean algebras given by $M_1,\dots ,M_i$.  All we need to do is show
that the matching on each $B_i \setminus (B_1 \cup \cdots \cup B_{i-1})$
is acyclic.  But if there were a cycle, let $\mu_i$ be the highest
precedence label to be inserted as a matching step in the cycle.
This would necessitate a downward step in the cycle shifting $\mu_i$ 
back into the interior of $I(M_r)$, but this would be preceded and
followed by matching steps inserting labels of lower precedence than
$\mu_i$.  This contradicts our greedy matching procedure, because it
would instead make the downward edge a matching edge inserting $\mu_i$,
since this has higher precedence than the matching step of either 
endpoint.  Thus, there are no cycles.
\EOP

\section{Rationality of Morse number
generating function}\label{survive-deg-d}

In this section we describe a finite state automaton that generates
exactly the language of label sequences for surviving critical cells, in
the case of a Gr\"obner basis of
degree $d$.  The existence of such a generating function again implies
the language is regular, and hence that 
the generating function for Morse numbers is a rational function
which gives upper bounds on the terms in the Poincare'-Betti series.
In contrast, for $d\ge 3$ the 
Poincare'-Betti series is not always rational.  The generating function
for Morse numbers does come close enough
to the Poincare'-Betti series to
achieve the vanishing of Betti numbers described by
Theorem ~\ref{deg-d}.
Due to the similarity of the finite state automaton
to the one given in the quadratic Gr\"obner basis case, less 
detail is provided here than in Section ~\ref{section-quad-applic}.

The states in the automaton keep track of the set of
previously encountered ILTs
and individual labels, in their order of most recent appearance.
Reading label sequences from top to bottom, 
the following are
the legal transitions from one state to another.  
\begin{enumerate}
\item
a single label $\lambda $
that is larger than its predecessor, i.e. which forms a 
descent with the label above it.  For the transition to be to a final
state, we require
the further property that $\lambda $ is separated 
from
each previously encountered ILT $I(a_1,\dots ,a_{d'})$ which satisfies
$a_1 < \lambda < a_{d'}$ and $\lambda a_2\cdots a_{d'},\lambda a_1\cdots
a_{d'-1}\not\in \init (I_{\Lambda })$.
by a label with which 
$\lambda $ does not commute.
\item
a single label which forms an ILT together with its predecessor, exactly
as in the $d=2$ case
\item
a collection $\{ a_1,\dots ,a_{d'-1}\}$  of labels with $d'>2$, which
together with the most recently encountered label $a_{d'}$ form an 
ILT $I(a_1,\dots ,a_{d'})$ such that (a) the 
labels $\{ a_2,\dots ,a_{d'}\} $ 
cannot all simultaneously shift upward into the interior of
higher ILTs to yield a critical cell which does not have any of the 
labels $\{ a_2,\dots ,a_{d'}\}$ individually as members of any 
non-essential set, and (b) no label above $I(a_1,\dots
,a_{d'})$ may shift downward into $I(a_1,\dots ,a_{d'})$ by a gradient 
path to yield a critical cell.  Such a transition leads to a final state.
\item
a collection of labels that collectively complete 
an ILT, with an allowable collection of interspersed
labels.  Allowable collections are those that arise as a result of 
concatenating ILTs, as described below, and those which may be within
the ILT in a surviving critical cell when there are multiple ILTs on
the same syzygy interval.  In this case the 
transition is to a non-final state, and we will justify below that there
are only a finite number of these transitions.
\end{enumerate}

The point is to use 
non-final states for label sequences for critical cells that are 
cancelled, if the concatenation of additional labels may yield a 
critical cell that is not cancelled.
to be cancelled.  The ``concatenating''
ILTs mentioned in the fourth type of transition come from situations
such as the following example.
\begin{examp}
Consider a label sequence $abcde$ where $abd$ and $cde$ are each 
Gr\"obner basis
leading terms.  The label $c$ cannot be shifted out of the ILT
$I(a,b,d)$ to yield a critical cell, because the ascend $(d,e)$ would 
no longer be part of a minimal skipped interval, so $c\not\in
S(a,b,d)$ for the label sequence $abcde$, though it would belong to 
$S(a,b,d)$ in the label sequence $abcd$. 
\end{examp}
Specifically, a pair of ILTs $M_1 = a_1\cdots a_r$ and 
$M_2 = b_1\cdots b_s$ are 
{\it concatenating} if either (a) $a_r$ divides $M_2$ with $a_r\ne b_s$ and
$a_1 < b_1$, or (b) $b_1$ divides $M_1$ with $b_1\ne a_1$ and 
$a_r < b_s$.

The fourth type of transition also accommodates the
matching procedure of 
Lemma ~\ref{mult-ilt}.

\begin{prop}
The automaton has finitely many states and transitions.
\end{prop}
\proof
There is a finite list of possible ILTs, even when we consider all 
possible label interspersions that could still allow the cell not to be
cancelled, i.e. from concatenating ILTs and from multiple ILTs on a 
single syzygy interval.  
This follows from the fact that the semi-group ring is 
finitely generated, and that each Gr\"obner basis leading term has finite
degree, so labels occurring in the interior of an ILT with multiplicity
greater than the Gr\"obner basis degree will always allow critical cell
cancellation.  The transitions out of a state are limited by the finite
list of labels.
\EOP

\begin{prop}
Word length
equals critical cell dimension shifted by two.
\end{prop}

\proof
Any label sequence which has more $I$-intervals than $J$-intervals will
be cancelled, unless there are two concatenating ILTs such that their
concatenation contains another ILT, causing three or more overlapping 
$I$-intervals in which one is discarded, in such a way that no labels
may be shifted from the interior of any of these ILTs without making the
cell non-critical.  
But in the case of this type of concatenation, where
two ILTs share labels and cover a third ILT, this means we can use just
the labels in these two ILTs for labelling the transitions in the 
finite state automaton, so we get the correct word length.
\EOP

Using the observations and propositions above, it is not hard to 
generalize the automaton from the $d=2$ case to obtain:

\begin{thm}
The surviving critical cells are labeled by 
the words of a regular language, with word length 
measuring cell dimension, shifted by two.  Thus, the generating
function for Morse numbers is a rational generating function which is 
determined by the given finite state automaton.  
\end{thm}

\section{Some remarks and open questions}\label{remarks-questions}

\begin{rk}
\rm{
When a variable does not appear in any syzygies, then it may be 
``factored out'' before starting our analysis, similar to the situation with
computing Tor groups directly.  
Specifically, if some $z_i$ does not appear in any generators 
of the toric ideal $I_{\Lambda }$ for 
$k[z_1,\dots ,z_n]/I_{\Lambda }$, then the partial order $\Lambda $
is the product of an infinite chain together with the 
poset of monomials ordered by divisibility in $k[z_1,\dots ,\hat{z_i},\dots ,
z_n]/I_{\Lambda }$.  Thus, any finite interval is the product of a 
finite chain together with a monoid poset interval $(0,\lambda )$ for 
the ring
$k[z_1,\dots ,\hat{z_i},\dots ,z_n]/I_{\Lambda }$; the order complex of 
such an interval is the suspension of the join of the order complexes 
for the two terms in
the product, so the suspension of the 
join of a simplex (i.e. the order complex of a chain)
with the order complex $\Delta (0,\lambda )$.
}
\end{rk}

\begin{rk}\label{non-unique-remark}
\rm{
We sometimes have gradient paths which reverse a decreasing sequence of
labels of length $d>2$ to produce an ILT, in which case the 
permutation on labels is not 321-avoiding.  We have not matched and
cancelled any such pairs of critical cells.
Theorem ~\ref{red-exp} shows there are at most two
gradient paths between 
a pair of critical cells related by such a reversal for lexicographic
discrete Morse functions; the proof of Theorem ~\ref{red-exp} generalizes
to those facet orders which satisfy the crossing condition, so in 
particular to content-lex facet orders.  
}
\end{rk}

\begin{thm}
Suppose $I_{\Lambda }$ has degree at most three.  Then
for each critical cell $\tau $ in our complex $\Delta^M$
of critical cells after cancellation, $\partial (\tau )$ is a linear
combination of critical cells of content strictly earlier than
$\tau $.
\end{thm}

\proof
Suppose $\tau , \sigma$, $\dim(\tau) = \dim(\sigma)+$, are surviving critical cells
with equal content and there is a gradient path from $\tau $ to $\sigma $.
Then $\tau ,\sigma $ each have no syzygy intervals with 
non-empty non-essential set.
Any gradient path from $\tau $ to $\sigma $ must
sort labels, but in such a way that $\sigma $ still has no syzygy intervals
with non-empty non-essential set.  This can only be accomplished by
reversing three or more descending labels to form a new
ILT.  This ILT must come from a Gr\"obner basis leading term of degree
exactly three, since pairs of labels comprising
degree 2 leading terms cannot be swapped without passing to an earlier
content class.  Lemma ~\ref{mult-ilt-deg3} ensures that 
the three or more labels must occur in a single string of 
descending labels within $\tau $, to avoid $\tau $ being cancelled by
virtue of a syzygy
interval with non-empty non-essential set.  

The ILT to be created cannot come from a Gr\"obner basis leading term
of degree greater than three, both because of the assumptions of our
theorem, and also because this would decrease critical cell dimension 
by more than one, implying $\sigma $ could not be in the image of 
the boundary map applied to $\tau $.  
Theorem ~\ref{red-exp} shows there are at most two gradient paths 
reversing three labels, resulting from the Coxeter relation $s_i s_{i+1}
s_i = s_{i+1} s_i s_{i+1}$ being applied at the conclusion of a reduced
expression.  But one may easily check that one will indeed get two 
gradient paths when we reverse three labels as required for $\sigma $ 
in the boundary of $\tau $, and furthermore, that these will be oriented
so that the two ways in which $\sigma $ is incident to $\tau $ will 
cancel.  Thus, $\sigma $ will appear with coefficient 0 in the boundary
of $\tau $.
\EOP

\begin{rk}
The following example shows that the Morse function
bound on which Tor groups vanish is sharp.
Consider 
$$k[\Lambda ] =
k[z_1,\dots ,z_{2d}]/(z_1\cdots z_d - z_{d+1}\cdots 
z_{2d}),$$ or
equivalently,  
$$k[z_1z_2,z_3z_4,\dots ,z_{2d-1}z_{2d},z_1z_{d+1},z_2z_{d+2},
\dots ,z_dz_{2d}].$$  This clearly has a Gr\"obner basis of degree $d$
and none of lower degree.  The interval $(1,z_1\cdots z_d)$ in 
$\Lambda $ is disconnected.
\end{rk}

\begin{qn}
Is there a nice description of the gradient paths between surviving
critical cells?  This would be needed for a completely explicit
description of the boundary maps in our resolution, since these
are sums over such gradient paths.
\end{qn}

\begin{qn}\label{veronese-qn}
Is it possible to improve our 
discrete Morse function into one that would provide
a combinatorial proof of the following theorem?
If an affine semi-group ring is standard graded, and its toric
ideal of syzygies has a Gr\"obner basis of degree $d$, 
then its $(d-1)$-st Veronese is Koszul.  
\end{qn}

In our setting, the above is equivalent to the
rank-selected subposet of $\Lambda $ consisting of exactly
the ranks divisible by $d-1$ being a Cohen-Macaulay poset.
Example ~\ref{examp-veronese} discusses the one situation in
which our critical cells skip more than $2d-3$ consecutive elements
of a saturated chain; this seems to be the main issue one 
would need to address to provide an affirmative answer, though one 
would also need to better understand the relationship between 
lexicographic discrete Morse functions and rank-selection or else to 
modify the Morse function to one for the rank-selected subposet.  
By Lemma ~\ref{mult-ilt-deg3}, the 
issue of skipping more than $2d-3$ consecutive ranks
does not arise for $d\le 3$.  

\begin{examp}\label{examp-veronese}
\rm{
There is only one
circumstance in which critical cells could skip more than $2d-3$
consecutive elements of a saturated chain, and this only may happen
in the $d>3$ case.  Namely, if there are 
distinct Gr\"obner basis
leading terms with the same initial and final labels, this may result
in overlapping Boolean algebras of critical cells, with cells with 
large syzygy intervals not necessarily cancelled.
}
\end{examp}

\begin{qn}
In [HRW], $\Tor $ groups related to
quotients of affine semi-group rings by monomial ideals are 
translated to homology of certain
relative complexes $\Delta (\lambda ,A)$, where $\lambda $ specifies a 
monoid poset interval and $A$ is a graphic subspace arrangement.  Does
our Morse function translate to this setting to provide useful
new information?
\end{qn}

\section*{Acknowledgments}

The authors thank Phil Hanlon, Mel Hochster, Ezra Miller 
and Vic Reiner for helpful conversations.

\end{document}